\documentclass{article}%
\usepackage{amsmath}
\usepackage{amsfonts}
\usepackage{amssymb}
\usepackage{graphicx}%
\providecommand{\U}[1]{\protect\rule{.1in}{.1in}}
\newtheorem{theorem}{Theorem}[section]

\newtheorem{conjecture}[theorem]{Conjecture}
\newtheorem{corollary}[theorem]{Corollary}

\newtheorem{definition}[theorem]{Definition}
\newtheorem{example}[theorem]{Example}

\newtheorem{lemma}[theorem]{Lemma}

\newtheorem{proposition}[theorem]{Proposition}

\begin{document}

\title{Strong approximation methods in group theory \\ LMS/EPSRC Short course}
\author{Nikolay Nikolov}
\date{Oxford, 10-14 September 2007}
\maketitle

\section{Introduction}

This course is concerned with linear groups $\Gamma\leq\mathrm{GL}_{n}(k)$
where $k$ is some field (usually of characteristic $0$). Linearity is one of
the most effective and well studied conditions one can put on a general
infinite group. The following are two of the most often used consequences of linearity.

\begin{enumerate}
\item[(i)] a finitely generated linear group $\Gamma$ is residually finite, and

\item[(ii)] if in addition char $k=0$, then $\Gamma$ is virtually torsion free.
\end{enumerate}

The first of these means that $\Gamma$ has many finite images, and one way to
study $\Gamma$ is to investigate these images (equivalently, the profinite
completion $\widehat{\Gamma}$ of $\Gamma$). One of the main objectives of this
course is the `Lubotzky alternative' for linear groups:

\begin{theorem}
\label{lubalt} Let $\Delta\leq\mathrm{GL}_{n}(k)$ be a finitely generated
linear group over a field $k$ of characteristic 0. Then one of the following holds:

\begin{enumerate}
\item[\emph{(a)}] the group $\Delta$ is virtually soluble, or

\item[\emph{(b)}] there exist a connected simply connected $\mathbb{Q}$-simple
algebraic group $G$, a finite set of primes $S$ such that $\Gamma
=G(\mathbb{Z}_{S})$ is infinite, and a subgroup $\Delta_{1}$ of finite index
in $\Delta$ such that every congruence image of $\Gamma$ appears as a quotient
of $\Delta_{1}$.
\end{enumerate}
\end{theorem}

Here $\mathbb{Z}_{S}=\mathbb{Z}[1/p\ |\ p\in S]$. In case (b) we can deduce
from the \emph{Strong Approximation Theorem} that $\Delta_{1}$ has many finite
images, in particular the groups $\prod_{i=1}^{k}G(\mathbb{F}_{p_{i}}) $
whenever $p_{1},\ldots,p_{k}$ are distinct primes outside $S$. Now, for all
but finitely many primes $p$ the group $G(\mathbb{F}_{p})$ is semisimple, in
fact it is a perfect central extension of a product of simple groups (of fixed
Lie type over $\mathbb{F}_{p}$). The simple groups of Lie type are very well
understood and this enables us to deduce properties of the profinite
completion $\widehat{\Delta}$ of $\Delta$.

For example, if $\Delta$ has polynomial subgroup growth then one can deduce
that case (b) of Theorem \ref{lubalt} is impossible and hence that $\Delta$ is
virtually soluble. Some more applications of Theorem \ref{lubalt} are given in
section \ref{applic} below. \medskip

In turn, when $\Delta$ is virtually soluble we have the \textbf{Lie-Kolchin
theorem}:

\begin{theorem}
\label{lk} Suppose that $\Delta\leq\mathrm{GL}_{n}(K)$ is a virtually soluble
linear group over an algebraically closed field $K$. Then $\Delta$ has a
triangularizable subgroup $\Delta_{1}$ of finite index; i.e. $\Delta_{1}$ is
conjugate to a subgroup of the upper triangular matrices in $\mathrm{GL}%
_{n}(K)$.
\end{theorem}

In fact if $\mathrm{char}K=0$ the index of $\Delta_{1}$ in $\Delta$ can be
bounded by a function of $n$ only (a theorem of Mal'cev and Platonov). This
has the corollary:

\begin{lemma}
\label{plat} Suppose that $\Delta$ is a finitely generated group which is
residually in the class of virtually soluble linear groups of degree $n$ in
characteristic $0$. Then $\Delta$ itself is virtually soluble.
\end{lemma}

We shall use this Lemma in the proof of Theorem \ref{lubalt}.

A common feature in the proofs of all these results is to take the
\emph{Zariski closure} $G=\overline{\Delta}$ of $\Delta$ in $\mathrm{GL}_{n}(K)$.
 This is a linear algebraic group and we can apply results from
algebraic geometry, number theory and the theory of arithmetic groups to study
$G$ and its dense subgroup $\Delta$. \medskip

The main object of this course is to understand the terminology appearing
above and develop the methods by which Theorems \ref{lubalt} and \ref{lk} can
be proved. These methods are useful in a variety of other situations involving
linear groups.

\section{Algebraic groups}

Throughout, $K$ will denote an algebraically closed field.

\subsection{The Zariski topology on $K^{n}$.}

A good reference for the material of this section (with proofs) is the book
\cite{AM} by Atiyah and Macdonald. For a brief introduction see also the
chapter `Linear algebraic Groups' in \cite{CMS}.

Let $K^{n}$ be the $n$-dimensional vector space over $K$. Given a subset $S
$ of the polynomial ring
\[
R:=K[x_{1},\ldots,x_{n}]
\]
define
\[
V(S)=\{x\in K^n \ |\ f(\mathbf{x})=0\quad\forall f\in S\}
\]
to be set of common zeros of $S$ in $K^n$.

It is easy to that $V(I)=V(S)$ for the ideal $I$ generated by $S$, that
\[
V(I)\cup V(J)=V(IJ)
\]
for ideals $I$ and $J$ of $R$, and that%
\[
\bigcap_{I\in\mathcal{F}}V(I)=V(\sum_{I\in\mathcal{F}}I)
\]
for any family $\mathcal{F}$ of ideals of $R$.

The Hilbert Basis Theorem says that each ideal $I$ of $R$ is finitely
generated and so each $V(S)$ can in fact be defined by finitely many
polynomial equations.

\begin{definition}
The \emph{Zariski topology} of $K^n$ has as its closed sets the sets $V(I)$ for
all ideals $I$ of $R$.
\end{definition}

The first basic result about the Zariski topology is the following 

\begin{proposition}[Exercise 3 on p. \pageref{exerc} ]\label{compact} The space $K^n$ with the Zariski topology is a compact topological space, in
fact it satisfies the descending chain condition on closed subsets. 
\end{proposition}

Note that the closed sets of $K$ coincide with its finite subsets 
(since a polynomial in one variable can only have finitely many roots). More generally 
the Zariski topology of $K^n$ is never Hausdorff, thus even though the space $K^n$ is compact 
one should be careful when 
applying familiar results from Hausdorff spaces.
\bigskip

\textbf{Example:} Let $V$ be the hyperbola given by the equation $x_1x_2=1$ in $K^2$.
Then $V$ is a closed, hence compact subset of $K^2$ but its projection on the $x_1$ axis is $K \backslash \{0\}$ which 
is not closed. So in the Zariski topology continuous images of compact sets are not always closed.  
\bigskip

The subsets $V(I)\subseteq K^{n}$ (with the subspace topology induced from the
Zariski topology on $K^{n}$) are called \emph{affine (algebraic) varieties}.
If $W$ is an affine variety, the \emph{coordinate ring} $R(W)$ of $W$ is the
algebra $R/J(W)$, where $J(W)$ is the ideal of $R$ consisting of all
polynomials vanishing on $W$. The ascending chain condition on ideals of $R$
(Hilbert's Basis Theorem) implies the descending chain condition (minimal
condition) for closed sets in $K^n$.

\begin{theorem}
\emph{(Hilbert's Nullstellensatz)}  $V(I)=\emptyset$ if and only if $I=R$.
\end{theorem}

In fact a more general result holds (see \cite{AM}, Chapter 7): if $W=V(I)$ is
an affine variety then $J(W)/I$ is the \emph{radical} of $I$, i.e.
\[
J(W)=\left\{  x\in R\ |\ x^{n}\in I\text{ for some }n\in\mathbb{N}\right\}  .
\]

The coordinate ring $R(W)$ can be considered as the set of morphisms of $W$
into the one-dimensional variety $K$. In general, a \emph{morphism} $F$ from
$W_{1}\subseteq K^{n_{1}}$ into $W_{2}\subseteq K^{n_{2}}$ is an $n_{2}$-tuple
$(f_{1},\ldots,f_{n_{2}})\in K[x_{1},\ldots,x_{n_{1}}]^{n_{2}}$ of polynomial
maps such that $F(W_{1})\subseteq W_{2}$. Any such morphism induces a
$K$-algebra homomorphism $F^{\ast}:R(W_{2})\rightarrow R(W_{1})$ defined by
$f\mapsto f\circ F$. Conversely, from the Nullstellensatz it can be shown that
every algebra homomorphism $F^{\ast}$ between $R(V_{2})$ and $R(V_{1})$ arises
in this way from a morphism $F:V_{1}\rightarrow V_{2}$. In this way the
category of affine varieties is anti-equivalent to the category of
reduced\footnote{$R$ is \emph{reduced} if it contains no non-zero nilpotent
elements.} finitely generated algebras over the algebraically closed field $K$.

\begin{definition}
A variety is \emph{irreducible} if it is not the union of two proper closed subsets.
\end{definition}

Since $V$ satisfies the minimal condition on closed subsets we can write every
variety $W$ as
\[
W=W_{1}\cup W_{2}\cup\cdots\cup W_{k},
\]
a union of irreducible varieties $W_{i}$. If we assume that the above
decomposition is irredundant, i.e. $W_{i}\nsubseteq W_{j}$ whenever $i\neq j
$, then it is in fact unique up to reordering of the $W_{i}$. These are then
called the \emph{irreducible components} of $W$.

For example if $W$ is the variety defined by the single equation
\[
x_{1}x_{2}(x_{1}x_{2}^{2}-1)=0
\]
then its irreducible components are the two lines with equations $x_{1}=0$,
$x_{2}=0$ and the curve defined by $x_{1}x_{2}^{2}=1$. \medskip

It is easy to see that a variety $W$ is irreducible if and only if $J(W)$ is a
prime ideal of $W$, i.e. if and only if the coordinate ring $R/J(W)$ is an
integral domain.

\begin{definition}
The \emph{dimension}, $\dim W$ of an irreducible variety $W$ is the Krull
dimension of $R(W)$. This is the transcendence degree of $R(W)$ over $K$, or
equivalently the maximal length $d$ of a chain of distinct nonzero prime
ideals $0\subset P_{1}\subset\cdots\subset P_{d}\subset R(W)$. The dimension
of a general affine variety is the maximal dimension of its irreducible components.
\end{definition}

As a consequence, a closed proper subset of an irreducible variety $W$ has
strictly smaller dimension than $W$.

\subsection{Linear algebraic groups as closed subgroups of $\mathrm{GL}_{n}(K)$.}

We identify the $n\times n$ matrix ring $M_{n}(K)$ with $K^{n^{2}}$, and use
$x_{ij}$ ($i,j=1,\ldots,n)$ as coordinates. Then the subgroup $\mathrm{SL}%
_{n}(K)$ of matrices with determinant $1$ forms an affine variety, defined by
the equation $\det(x_{ij})=1.$

\begin{definition}
A \emph{linear algebraic group} over $K$ is a Zariski-closed subgroup of
$\mathrm{SL}_{n}(K)$ for some $n$.
\end{definition}

\textbf{Notes}:

\begin{enumerate}
\item The two maps $(x,y)\mapsto xy$ and $x\mapsto x^{-1}$ from $G\times G$
(resp. $G$) to $G$ are morphisms of affine varieties.

\item There are more general algebraic groups which are not linear. In this
course we shall be concerned only with linear algebraic groups and `algebraic
group' will always mean `linear algebraic group'.

\item The definition we have given is different from the standard one but
equivalent to it: one usually defines a linear algebraic group to be an affine
variety with maps of group multiplication and inverses which are morphisms of
varieties. It can be shown that every such group is in fact isomorphic to a
closed subset of some $\mathrm{SL}_{n}(K)$. See the `Linear algebraic groups'
chapter in \cite{CMS}.
\end{enumerate}

\bigskip

A \emph{homomorphism} of linear algebraic groups $f:G\rightarrow H$ is a group
homomorphism which is also a morphism between varieties. i.e. $f$ is given by
polynomial maps on the realizations of $G\subset M_{n_{1}}(K)$ and $H\subset
M_{n_{2}}(K)$.

The group $\mathrm{GL}_{n}(K)$ is isomorphic to a closed subgroup of
$\mathrm{SL}_{m+1}(K),$ by the mapping%
\[
g\mapsto\left(
\begin{array}
[c]{cc}%
g & 0\\
0 & (\det g)^{-1}%
\end{array}
\right)  .
\]
In this way we consider $\mathrm{GL}_{n}(K)$ as a linear algebraic group. It
is clear that every linear algebraic group is isomorphic to a closed subgroup
of $\mathrm{GL}_{n}(K)$ for some $n$.

\subsubsection{Basic examples}

For an integer $n\geq2$ consider the following subgroups of $\mathrm{SL}_{n}(K)$:

\begin{itemize}
\item The group of (upper or lower) unitriangular matrices,

\item The upper (upper or lower) triangular matrices,

\item The diagonal matrices, or more generally

\item The monomial matrices.
\end{itemize}

It is clear that these are closed subgroups of $\mathrm{SL}_{n}(K)$ and so are
algebraic groups.

Note that when $n=2$ the first example is isomorphic to the additive group of
the field $K$, while the third one is isomorphic to the multiplicative group
of $K$. In this way $(K,+)$ and $(K,\times)$ become linear algebraic groups.
The first one is denoted by $\mathbb{G}_{+}$ and the second by $\mathbb{G}%
_{\times}$. In can be shown that these are the only connected algebraic groups
of dimension $1$.

Another family of examples arise from linear groups preserving some form. For
example if $(\mathbf{u},\mathbf{v})=\mathbf{u}^{T}P\mathbf{v}$ is a bilinear
form on the vector space $V=K^{n}$, then the group 
$G\leq\mathrm{GL}(V)$ preserving $(-,-)$ can be described as those matrices $X$ in
$\mathrm{GL}_{n}(K)$ such that $X^{T}PX=P$. This is a collection of $n^{2}$
polynomial equations on the coefficients of $X=(x_{ij})$ and so $G$ is an
algebraic group. Examples are the symplectic group $\mathrm{Sp}_{n}(K)$ ($n$
even) and the special orthogonal group $\mathrm{SO}_{n}(K)$.

\subsubsection{Basic properties of Algebraic Groups}

\begin{theorem}
\emph{(see \cite{Hu}, chapter II)}  Let $f:G\rightarrow H$ be a homomorphism of algebraic
groups. Then

\begin{enumerate}
\item $\mathrm{Im}(f)$ is a closed subgroup of $H$ and $\ker(f)$ is a closed
subgroup of $G$.

\item $\dim G=\dim\ker(f)+\dim\mathrm{Im}(f)$.
\end{enumerate}
\end{theorem}

Recall that a topological space is \emph{connected} if cannot be written as
the disjoint union of two proper closed (equivalently open) subsets. Clearly
an irreducible variety is connected. It turns out that for algebraic groups
the converse is also true and so the two concepts coincide. To see this,
suppose that $G$ is a connected algebraic group. Let $G=V_{1}\cup\cdots\cup
V_{k}$ be the decomposition of $G$ into irreducible components. This
decomposition is unique up to the order of the $V_{i}$, therefore the action
of $G$ by left multiplication permutes the components $V_{i}$. Without loss of
generality suppose that $1\in V_{1}$. Let
\[
G_{1}=\mathrm{Stab}_{G}(V_{1}):=\{g\in G\ |\ gV_{1}=V_{1}\}.
\]

Clearly $G_{1}$ is a closed subgroup of finite index $k$ in $G$, so it is both
open and closed, as are each of its cosets in $G$. Since $G$ is connected we
must have $G=G_{1}$, so $k=1$ and $G$ is irreducible.

The above argument easily shows that more generally the connected component of
the identity $G^{\circ}$ of $G$ is a closed irreducible normal subgroup of
finite index in $G$; it may be characterized as the smallest closed subgroup
of finite index in $G$.

\begin{lemma}
\emph{(see \cite{Wef}, \textbf{14.15} or \cite{Hu}, \S 7.5)}  If $(H_{i})_{i\in I}$ is a family of
closed connected subgroups of $G$ then the subgroup $\langle H_{i}|\ i\in
I\rangle$ generated abstractly by the $H_{i}$ in $G$ is closed and connected.
\end{lemma}

In particular if $H_{1}$ and $H_{2}$ are two closed connected subgroups of $G$ such that
$H_{1}H_{2}=H_{2}H_{1}$ (e.g., if either of $H_{1}$ or $H_{2}$ is normal in
$G$) then $H_{1}H_{2}$ is a closed connected subgroup of $G$. In general
if $H_1, H_2 $ are closed subgroups with $H_1H_2=H_2H_1$ and having connected components
$H_1^0$, $H_2^0$ respectively, then $H_1H_2$ is
a finite union of closed sets $hH_1^0H_2^0h'$ for some $h,h' \in G$
and so is a closed subgroup of $G$.

\begin{theorem}
\emph{(Chevalley; see \cite{Hu}, chapter IV)}  If $H$ is a closed normal subgroup of $G$ then
the quotient $G/H$ can be given the structure of a linear algebraic group, so
that the quotient map $G\rightarrow G/H$ is a homomorphism of algebraic groups.
\end{theorem}

\subsubsection{Fields of definition and restriction of scalars.}

A variety $V(I)$ is said to be \emph{defined over} a subfield $k\subset K$ if
the ideal $I$ is generated (as an ideal of $R$) by polynomials with
coefficients in $k$. When the field $k$ is separable (which is always the case
if $k$ has characteristic $0)$ there is a useful criterion for $V$ to be
defined over $k$:

\begin{lemma}
Let $W=V(S)$ be a variety. For $\sigma\in\mathrm{Gal}(K/k)$ define the variety
$W^{\sigma}$ to be $V(S^{\sigma})$, i.e., the zero set of the ideal
$S^{\sigma}$ of $R$. Then $W$ is defined over $k$ if and only if $W=W^{\sigma
}$ for all $\sigma\in\mathrm{Gal}(K/k)$.
\end{lemma}

Similarly, a homomorphism $f:G\rightarrow H$ between two algebraic groups is
$k$-defined if all the coordinate maps defining $f$ are polynomials with
entries in $k$.

Now let $G\leq\mathrm{GL}_{n}(K)$ be an algebraic group and let $\mathcal{O}$
be a subring of $K$. The group of $\mathcal{O}$\emph{-rational points} of $G$
is defined to be $\mathrm{GL}_{n}(\mathcal{O})\cap G$ and is denoted by
$G_{\mathcal{O}}$.

Suppose that $G$ is defined over some subfield $k$ of $K$ which is a finite
extension of $k_{0}$. In this course we shall study the groups $G_{k}$ and
sometimes we prefer to reduce the situation to a smaller subfield $k_{0}$
(which will usually be $\mathbb{Q}$).

There is a standard procedure for doing this, called `restriction of scalars'.
This associates to $G$ another algebraic group $H\leq\mathrm{GL}_{nd}(K)$
where $d=(k:k_{0})$; here $H$ is defined over $k_{0}$ and satisfies $H_{k_{0}%
}=G_{k}$. The algebraic group $H$ is denoted $\mathcal{R}_{k/k_{0}}(G)$.
Before presenting the general construction let us study a simple special case
which illustrates the idea.

Suppose that $G$ is the multiplicative group of the field $(K,\times)$. This
is defined over the integers $\mathbb{Z}$ (i.e. it can be defined by
polynomials over $\mathbb{Z}$.) Let $k$ be a number field, i.e. a finite
extension field of $\mathbb{Q}$. The group $G_{k}$ is clearly the
multiplicative group $k^{\ast}$ of the field $k$. We want to find a
$\mathbb{Q}$-defined algebraic group $H$ such that its group $H_{\mathbb{Q}}$
of $\mathbb{Q}$-rational points is isomorphic (as an abstract group) to
$G_{k}$.

To find $H$ we identify $k$ with the vector space $\mathbb{Q}^{d}$ by choosing
a basis $a_{1},\ldots a_{d}$ for $k$ over $\mathbb{Q}$, and consider the
regular representation of $k$ acting on itself by left multiplication. This
gives an algebra monomorphism $\rho:k\rightarrow M_{d}(\mathbb{Q})$ and so
$\rho(k)$ is a $d$-dimensional subspace of $M_{n}(\mathbb{Q})$. This can be
defined as the zeroes of some $s=d^{2}-d$ linear functionals $F_{1}%
,\ldots,F_{s}:M_{n}(\mathbb{Q})\rightarrow\mathbb{Q}$. Therefore we can define
the algebraic variety $H$ as the set of zeros of $F_{1},\ldots,F_{s}$ in
$\mathrm{GL}_{d}(K)$. Then clearly $H_{\mathbb{Q}}=G_{k}$ and the only thing
that has to be checked is that $H$ is a group, i.e. the variety $H$ is closed
under matrix multiplication and inverses. This can be expressed as the
vanishing of certain polynomials in the coordinates $x_{ij}$. If one of these
polynomials is nontrivial it will be nontrivial for some rational values of
its arguments. But we certainly know that $H_{\mathbb{Q}}$ is closed under
multiplication and inverses since it is equal to the multiplicative group
$k^{\ast}$. So $H$ is indeed an algebraic group.

There is another way to view the algebraic group $H$ just constructed. Let
$\sigma_{1},\ldots,\sigma_{d}$ be the $d$ embeddings of $k$ into the
algebraically closed field $K$. For an element $h=\sum_{i=1}^{d}%
x_{i}a_{i}\in k$ with $x_{i}\in\mathbb{Q}$ consider
\[
\lambda(h)=(\lambda_{1}(h),\ldots,\lambda_{d}(h)),
\]
where
\[
\lambda_{j}(h)=\sum_{i=1}^{d}x_{i}\sigma_{j}(a_{i})=\sigma_{j}(h).
\]
The condition $\det\rho(h)\not =0$ is equivalent to $\prod_{j}\lambda
_{j}(h)\not =0$. If $k=\mathbb Q (\alpha_1)$ where 
$\alpha_1$ has minimal polynomial $p(x)=(x-\alpha_1) \cdots (x-\alpha_d)$
over $\mathbb Q$ then $k \simeq \mathbb Q[x]/(p(x))$. We can extend $\lambda$ from $k$ to $k \otimes \mathbb Q \mathbb C$ and then  
\begin{equation} \label{chin} k \otimes_\mathbb Q \mathbb C \simeq \frac{ \mathbb C [x]}{(p(x))} 
\simeq \bigoplus_{i=1}^d \frac{\mathbb C[x]}{(x-\alpha_i)},
\end{equation} where the second
isomorphism comes from the Chinese remainder theorem and coincides with $\lambda$. Thus $\lambda \circ \rho^{-1}$ provides a $K$-isomorphism of $H$ with the direct product
$(\mathbb{G}_{\times})^{d}$ of $d$ copies of the multiplicative group
$\mathbb{G}_{\times}$. \label{chinese}

\bigskip

In general we are given a $k$-defined algebraic group $G\leq\mathrm{GL}_{n}(K)$.
Consider again an embedding $\rho:k\rightarrow M_{d}(k_{0})$ given
by the regular representation of $k$ acting on itself. Again the subspace
$\rho(k)\subset M_{d}(k_{0})$ is defined by some set of say $r$ linear
equations $F_{i}(y_{ab})$ in the matrix entries $y_{ab}$ ($1\leq a,b\leq d$
and $1\leq i\leq r$). Suppose that $G$ was defined as a variety by the $l$
polynomials $P_{j}(z^{st})$ in the entries of the matrix 
$(z^{st})\in M_{n}(K)$ ($j=1,\ldots,l,\ \ 1\leq s,t\leq n$).

Now the algebraic group $H=\mathcal{R}_{k/k_{0}}(G)$ is defined by the
following two families equations in the $(nd)^{2}$ variables $z_{ab}^{st}$:

The first family is
\[
P_{j}((z_{ab}^{st})_{a,b})=0\in M_{d}(K),\quad j=1,2,\ldots,l,
\]
i.e., we replace each variable $z^{st}$ in the original polynomial $P_{j}$
with a matrix $(z_{ab}^{st})_{a,b}\in M_{d}(K)$. Note that each $P_{j}$ gives
$d^{2}$ polynomial equations in $K$, one for each entry of the matrix in
$M_{d}(K)$.

The second family is
\[
F_{i}((z_{ab}^{st})_{a,b})=0,\ \ i=1,\ldots,r
\]
for each pair $(s,t)$ with $1\leq s,t\leq n$.

A typical example is the group
\[
G=\left\{
\begin{pmatrix}
a & 2b\\
b & a
\end{pmatrix}
\ |\quad a^{2}-2b^{2}\not =0\right\}
\]
which is the restriction of scalars $\mathcal{R}_{\mathbb{Q}(\sqrt
{2})/\mathbb{Q}}\mathbb{G}_{\times}$. Here, $G$ is $K$-isomorphic to
$\mathbb{G}_{\times}\times\mathbb{G}_{\times}$ via the map $%
\begin{pmatrix}
a & 2b\\
b & a
\end{pmatrix}
\mapsto(a+ib,a-ib)$, but this isomorphism is not $\mathbb{Q}$-defined.
\bigskip

It is easy to see that if we have a $k$-defined morphism $f:G\rightarrow T$
between two $k$-defined linear algebraic groups then this induces a $k_{0}%
$-defined morphism
\[
\mathcal{R}_{k/k_{0}}(f):\mathcal{R}_{k/k_{0}}(G)\rightarrow\mathcal{R}%
_{k/k_{0}}(T).
\]
In this way $\mathcal{R}_{k/k_{0}}$ becomes a functor between the category of
$k$-defined groups and morphisms and $k_{0}$-defined groups and morphisms.

\subsubsection{The Lie algebra of $G$}

There is a standard way to associate a Lie algebra $L(G)$ to any connected
linear algebraic group $G$, so that the map $L:G\mapsto L(G)$ is an
equivalence of categories. More precisely the following holds (see III of
\cite{Hu}):

\begin{itemize}
\item If $f:G\rightarrow H$ is a homomorphism of algebraic groups then there
is a uniquely specified homomorphism $L(f):L(G)\rightarrow L(H)$ between their
Lie algebras.

\item In particular, for any given $g\in G$ the map $x\mapsto g^{-1}xg$ is an
automorphism of $G$ and this gives rise to a Lie algebra automorphism denoted
$\mathrm{Ad}g:L(G)\rightarrow L(G)$. In this way we get a homomorphism of
algebraic groups $\mathrm{\ Ad:}G\rightarrow\mathrm{Aut}L(G)$, and it is easy
to see that $\ker\mathrm{Ad}=\mathrm{Z}(G)$, the centre of $G$.

\item If $H$ is a closed (normal) subgroup of $G$ then $L(H)$ is a Lie
subalgebra ( resp. an ideal) of $L(G)$.

\item If $G$ is defined over a subfield $k$ of $K$ then $L(G)$ is also defined
over $k$, i.e., it has a basis such that the structure constants of the lie
bracket multiplication are elements of $k$. Moreover if the morphism $f: G
\rightarrow H$ is $k$-defined then so is the Lie algebra homomorphism $L(f)$.

\item The dimension of $L(G)$ (as a vector space over $K$) is equal to the
dimension of the algebraic group $G$.
\end{itemize}

In general if $G$ is not connected we define $L(G)$ to be equal to $L(G^{0})$
where $G^{0}$ is the connected component of the identity in $G$. \medskip

Now a linear algebraic group $G$ is an affine subset of $M_{n}(K)$ so it is
defined by an ideal $I \vartriangleleft R$ of the polynomial ring $K[X_{11},
\ldots, X_{nn}]$. In this setting there is a concrete description of $L(G)$.
It is a Lie subalgebra of the Lie algebra $M_{n}(K)$ with the Lie bracket
\[
[A,B]=AB-BA.
\]

As a vector space $L(G)$ is the tangent space at the identity element $e \in
G$. In our situation this is defined as follows.\bigskip

For a polynomial $P\in R=K[(x_{ij})]$ and $g=(g_{ij})\in G\leq M_{n}(K)$ let
$\partial P_{g}$ be the linear functional on $n^{2}$ variables $X_{ij}$
defined as follows
\[
\partial P_{g}:M_{n}(K)\rightarrow K,\quad\partial P_{g}((X_{ij})_{i,j}%
):=\sum_{i,j}\left(  \frac{\partial P}{\partial x_{ij}}(g_{ij})\cdot
X_{ij}\right)
\]

Then $L(G)$ is the subspace of $M_{n}(K)$ of common solutions to the
equations
\[
\partial P_{e}=0,\quad\forall P\in I,
\]
where $e=\mathrm{Id}_{n}$ is the identity matrix in $G\leq\mathrm{GL}_{n}(K)$.

In fact we don't need to check infinitely many equations. By the Hilbert basis
theorem the ideal $I$ is finitely generated, say by polynomials $P_{1}%
,\ldots,P_{t}$. Then $L(G)$ is the common zeroes of the linear functionals
$\partial(P_{i})_{e}=0\ (i=1,\ldots,t)$.

\subsubsection{Connection with Lie algebras of Lie groups}
Let $G \leq GL_n$ be a linear algebraic group and suppose 
that $k$ is a complete field, for example $\mathbb C, \mathbb R$ or the field of $p$-adic numbers 
$\mathbb Q_p$ (see example \ref{qp} below). We have another topology on $\mathrm{GL}_n(k)$ by considering it as a subset of the 
topological space $M_n(k)=k^{(n^2)}$. In this way the group $G_k$ of $k$-rational points is a topological group, 
by virtue of being a closed subgroup of $\mathrm{GL}_n(k)$. In fact $G_k$ is a complex or real Lie group when 
$k=\mathbb C$ or $\mathbb R$, and is a $p$-adic analytic group when $k= \mathbb Q_p$.
In this section we shall use the term \emph{analytic group} to refer to either a Lie group or a $p$-adic analytic group.
 
There is a standard way 
to associate a Lie algebra $L(\mathcal G)$ to any (complex or real) Lie group $\mathcal G$ and as 
explained in \cite{bennotes} such a Lie algebra exists for any $p$-adic analytic group. One uniform way to define them 
is as the tangent space at the identity of $\mathcal G$. The Lie bracket is the differential of the commutator map in $\mathcal G$.     

The following Proposition is thus almost self evident. 

\begin{proposition} \label{liep} If the field $k$ is one from $\mathbb C, \mathbb R$ or $\mathbb Q_p$ then the $k$-rational points of $L(G)$, 
(namely $L(G)_k=L(G) \cap M_n(k)$) coincide with the Lie algebra of the analytic group 
$\mathcal G= G_k$.
\end{proposition}

For later use we record another basic result. First observe that when we have an analytic group 
$\mathcal G \leq \mathrm{GL}_n(k)$ with a faithful linear representation in $\mathrm{GL}_n(k)$ then 
we can also consider the Zariski topology
on $\mathcal G$ as a subset of $\mathrm{GL}_n$.

\begin{proposition} \label{ideals} Suppose that the group $H$ is a Zariski dense subgroup of the analytic group 
$G_k \leq \mathrm{GL}_n(k)$ 
for $G$ and $k$ as above. Let $\mathrm{Ad}$ be the adjoint action of $G$ on its Lie algebra $L(G)$. Then 
$\mathrm{Ad}(H)_k$ and $\mathrm{Ad}(G)_k$ have the same span in the vector space 
$\mathrm{End}_k L(G)_k$ over $k$.

Moreover, when $H$ is an analytic Zariski-dense in $G$ the Lie algebra $L(\mathcal H)$ of $\mathcal H$ is an ideal of the Lie algebra 
$L(G)_k$ of $G_k$.
\end{proposition}

\textbf{Proof.}  The adjoint action of $G$ on $L(G)$ is given by a set of 
polynomials (it coincides with the conjugation action of $G$ on $L(G)$ as a subset of $M_n(K)$ and so the map 
$\mathrm{Ad}: G \rightarrow \mathrm{End}_K(L(G))$ is morphism of algebraic varieties, 
hence a representation of $G$ as an algebraic group. Since $H$ is Zariski-dense in $G_k$ it follows that
$\mathrm{Ad} (H)$ is Zariski-dense in $\mathrm{Ad}(G_k)$ as subsets of $\mathrm{End}_k L(G)_k$. Since a 
vector subspace of $\mathrm{End}_k L(G)_k$ is Zariski closed the first part of the Proposition follows immediately.
\medskip

By a standard result of Lie theory $L(\mathcal H)$ is a Lie subalgebra of 
$L(G)_k$ which is $\mathrm{Ad}(\mathcal H)$-invariant. Now the stabilizer 
$\mathrm{Stab}(L(\mathcal H))$ of $L(\mathcal H)$ in $\mathrm{End}_k L(G)_k$ is a subspace of 
$\mathrm{End}_k L(G)_k$. Since this stabilizer contains $\mathrm{Ad}(\mathcal H)$ it should also 
contain $\mathrm{Ad}(G_k)$ by the first part of the Proposition.  
Therefore $L(\mathcal H)$ is $\mathrm{Ad}(G_k)$ invariant and so it is an ideal of the Lie algebra $L(G)_k$ of $G_k$.

\bigskip

\textbf{Note:} The Lie algebra is a local tool, it was only defined from a neighbourhood of the identity of an analytic 
group $\mathcal G$. 
So it is the same for any open subgroup of $\mathcal G$. 
In particular any subgroup of finite index in $G_{\mathbb Z_p}$ is a compact open subgroup of $G_{\mathbb Q_p}$ and 
hence all of these groups share the same Lie algebra as analytic groups. 
In fact this property characterises the open subgroups of analytic groups:

\begin{proposition}\label{openlie} Suppose that the analytic group $\mathcal H$ is a closed subgroup of the analytic group 
$\mathcal G$. Then $\mathcal H$ is an open subgroup of $\mathcal G$ if an only if $\mathcal H$ and $\mathcal G$ have the same Lie algebra. 
In particular when $\mathcal G$ is compact this happens if and only if $\mathcal H$ has finite index in $\mathcal G$. 
\end{proposition} 

\bigskip

As in the theory of Lie groups the Lie algebra is a very useful tool in the
study of algebraic groups. This is best seen in the classification of the
simple algebraic groups in next section, but we can already give a nontrivial application.

\begin{proposition}
A connected linear algebraic group $G$ of dimension less than 3 is soluble. If
$\dim G = 1$ then $G$ is abelian.
\end{proposition}

Indeed $L=L(G)$ is a Lie algebra of dimension at most 2 as a vector space over
$K$ and it is easy to see that $L$ must be soluble. If $\dim L=1$ then $L$ is
abelian and then so is $G$.

Note that even at this small dimension we see that two connected groups (for
example $\mathbb{G}_{+}$ and $\mathbb{G}_{\times}$) may have the same Lie
algebra and still be non-isomorphic. However the simply connected semisimple
groups are indeed uniquely determined by their Lie algebras as we shall see in the following section.

\subsection{Semisimple algebraic groups. The classification of simply
connected algebraic groups over $K$}

\begin{definition}
A connected algebraic group is called \emph{semisimple} if it has no
nontrivial closed connected soluble normal subgroups.
\end{definition}

In general, an algebraic group $G$ has a unique maximal connected soluble
normal subgroup. This is called the (soluble) \emph{radical} and denoted
Rad$(G)$. The group $G/\mathrm{Rad}(G)$ is then semisimple.

\begin{definition}
A connected algebraic group is \emph{simple} if it is nonabelian and has no
proper nontrivial connected normal subgroups.
\end{definition}

This implies that every closed proper normal subgroup of $G$ is central and
finite (\emph{Exercise}: prove this!).

\begin{theorem}
\label{ss} A semisimple group $G$ is a central product
\[
G\cong S_{1}\circ S_{2}\circ\cdots\circ S_{l}%
\]
of simple algebraic groups $S_{i}$. The factors in this product are unique up
to reordering.
\end{theorem}

\noindent Recall that a central product $S_{1}\circ S_{2}\circ\cdots\circ
S_{l}$ is a quotient $P/N$ of the direct product $P=S_{1}\times\cdots\times
S_{l}$ by a central subgroup $N$ intersecting each $S_{i}$ trivially.

So in order to understand semisimple algebraic groups it is sufficient to
understand simple algebraic groups and their central extensions. \medskip

Analogous definitions apply relative to any field of definition $k$. A
connected nonabelian algebraic group defined over $k$ is $k$\emph{-simple}
(resp. $k$\emph{-semisimple}) if it has no nontrivial closed connected proper
(resp. soluble) normal subgroups defined over $k$. Again a $k$-semisimple
group is $k$-isomorphic to a central product of $k$-simple groups which are
unique up to reordering.

When we speak of simple/semisimple groups without indicating the field the
understanding is that it is $K$. In this case $G$ is called \emph{absolutely}
simple (resp. semisimple). \textbf{Warning}: a $k$-simple group need not be
absolutely simple (though it is semisimple); see Example \ref{quat} below.\bigskip

The classification of absolutely simple algebraic groups mirrors entirely the
classification of the finite dimensional simple Lie algebras over $K$. Indeed
a simple group $G$ has finite centre and $G/\mathrm{Z}(G)$ embeds via
$\mathrm{Ad}$ as a group of automorphisms of its Lie algebra $L=L(G)$. So once
the algebra $L(G)$ is known the group $G$ is determined up to an \emph{isogeny}
as a
closed subgroup of $\mathrm{Aut}(L)$. More precisely we have the following classification
theorem.

\begin{theorem}
\emph{(Chevalley)}\label{K-simple} For each Lie type
$\mathcal{X}$ from the list
\[
A_{n}\ (n\geq1),\ B_{n}\ (n\geq2),\ C_{n}\ (n\geq3),\ D_{n}\ (n\geq
4),\ G_{2},\ F_{4},\ E_{6},\ E_{7},\ E_{8}%
\]
there are two distinguished simple groups of type $\mathcal{X}$: the so-called
\emph{simply connected} group $G_{sc}$ and the \emph{adjoint group}
$G_{ad}=G_{sc}/\mathrm{Z}(G_{sc})$. Every simple group of type $\mathcal{X}$
is an image of $G_{sc}$ modulo a finite central subgroup $T$. Such a quotient
map $\pi:\ G\rightarrow G/T$ is called a (central) \emph{isogeny}; all the
groups of the same type $\mathcal{X}$ form one isogeny class.

Every simple algebraic group belongs to exactly one of the isogeny classes
described above.
\end{theorem}

The proof of the uniqueness of the isogeny classes can be found in \cite{Hu} Chapter XI (see 
\textbf{Theorem'} in 32.1 there). Their \emph{existence} is discussed briefly in \cite{Hu} 33.6 and the construction of the groups of adjoint types is given in \cite{carter}.  

Examples of simply connected groups are $\mathrm{SL}_{n}(K)$ of type $A_{n-1}$
and $\mathrm{Sp}_{2n}(K)$ (type $C_{n}$). The group $\mathrm{SO}_{n}(K)$ is
simple of type $B_{(n-1)/2}$ or $D_{n/2}$ (depending on whether $n$ is even or
odd) but is not simply connected: its universal cover (i.e. the simply
connected group in its isogeny class) is $\mathrm{Spin}_{n}(K)$, the so-called
spinor group.

We extend the definition of `simply connected' to the semisimple groups:

\begin{definition}
A semisimple group is simply connected if it is a direct product of simply
connected simple groups.
\end{definition}

From Theorem \ref{K-simple} it now follows that each semisimple group is the
image by a central isogeny of a unique simply connected semisimple
group.\bigskip

In general the $k$-simple algebraic groups are not so easy to describe. In the
first place the radical of such a group is defined over $k$ and so it must be
trivial. Therefore a $k$-simple (even a $k$-semisimple) group is also
absolutely semisimple.

The next example gives a $\mathbb{Q}$-simple group which is not absolutely
simple.

\begin{example}
\label{quat} 
Let $H=\mathcal{R}_{\mathbb{Q}(i)/\mathbb{Q}}\mathrm{SL}_2$ be the restriction of scalars
of $\mathrm{SL}_2$ (defined over $\mathbb Q$) from $\mathbb{Q}(i)$ to $\mathbb{Q}$.
Then $H$ is $\mathbb Q$-simple.
\end{example}

Indeed, by Exercise 6 on page 23 we  see that $H$ is $\mathbb Q (i)$-isomorphic to 
$\mathrm{SL}_{2} \times \mathrm{SL}_2$ via an isomorphism $\rho$, say. Composing $\rho$ with 
complex conjugation $\tau$ has the effect of swapping the two factors $\mathrm{SL}_2$. 
It follows that none of these two factors is $\mathbb Q$-defined as a subgroup of $H$. 
Now if $H$ had a proper $\mathbb Q$-defined normal subgroup $L$ then $L$ must coincide with 
one of the two factors $SL_2$ but they are not defined over $\mathbb Q$. Contradiction, 
hence $H$ is $\mathbb Q$-simple.

\bigskip

Suppose now that $G$ is a $k$-simple, connected and simply connected group.
This means that over $K$ our group $G$ is isomorphic to a direct product
$\prod_{i} H_{i}$ of $K$-simple simply connected group $H_{i}$. It happens
that each of $H_{i}$ is defined over some finite Galois extension $k_{1}$ of
$k$ and we have that $G$ is $k$-isomorphic to the restriction of scalars
$\mathcal{R}_{k_{1}/k} H$ where $H=H_{1}$, say.

The group $H$ is $K$-simple so over $K$ it is isomorphic to one of the (simply
connected) groups listed in Theorem \ref{K-simple} but we need to classify
such groups up to $k_{1}$-isomorphism.

For example the group
\[
\mathrm{SO}_{2}=\left\{
\begin{pmatrix}
a & b\\
-b & a
\end{pmatrix}
\ |\ \ a^{2}+b^{2}=1\right\}
\]
is isomorphic to the multiplicative group $\mathbb{C}_{\times}$ over $K$ but
this isomorphism is not defined over the real subfield $\mathbb{R}$. \medskip

The $k_{1}$-isomorphism classes of groups which are $K$-isomorphic to $H$ are
called the $k_{1}$-forms of $H$. These are classified by the non-commutative
1-cohomology set $H^{1}(\mathrm{Gal}(K/k_{1}),\mathrm{Aut}H_{ad})$. For
example the unitary group $\mathrm{SU}_{n}$ is isomorphic to $\mathrm{SL}_{n}$
over $K=\mathbb{C}$ but not over $\mathbb{R}$ and these are the only two real
forms of $\mathrm{SL}_{n}$. Similarly the group $G$ in Exercise 7 on page
\pageref{ex7} is a
$\mathbb{Q}(i)$-form of $\mathrm{SL}_{2}$. For more details we about the
classification the $\mathbb{Q}$-forms of classical groups we refer to
\cite{PR}, Chapter 2.

\subsection{Reductive groups}

A class of groups which is more general than semisimple groups but which
shares some of their nicer properties is the \emph{reductive groups}. For
example $\mathrm{GL}_{n}(K)$ is not semisimple but still a very important
group which we would like to include in our theory.

\begin{definition}
An element $g$ of a linear algebraic group $G\leq\mathrm{GL}_{n}(K)$ is called
\emph{unipotent} (resp. \emph{semisimple}) if $g$ is unipotent (resp.
diagonalizable) as a matrix in $\mathrm{GL}_{n}(K)$. This definition is
independent of the choice of the linear representation of $G$. The group $G$
is \emph{unipotent} if it consists of unipotent elements.
\end{definition}

For example $\mathbb{G}_{+}$ is unipotent.

Now it can be shown that an algebraic group $G$ has a unique maximal normal
unipotent subgroup. This is the \emph{unipotent radical} of $G$ and is denoted
$R_{u}(G)$. The group $G$ is called \emph{reductive} if $R_{u}(G)=1$. One
obvious example of a reductive group is a torus.

\begin{definition}
An algebraic group $T$ is a \emph{torus} if it is isomorphic to a direct product
$\mathbb{G}_{\times}^{m}$. The \emph{rank} of $T$ is the number $m$ of
direct factors $\mathbb{G}_{\times}$. The torus $T$ is called \emph{$k$-split}
if there is a $k$-defined isomorphism $T \rightarrow \mathbb{G}_{\times}^{m}$  
\end{definition}

The group of diagonal matrices in $\mathrm{GL}_{n}(K)$ is a torus of rank $n$.

\begin{theorem}
\label{reductive} A connected reductive group $G$ is a product $G=TS$ of a
torus $T$ and a semisimple subgroup $S$ such that $[T,S]=1$ and $T\cap S$ is
finite. The subgroups $T$ and $S$ are unique.
\end{theorem}

For example $\mathrm{GL}_{n}(K)$ is reductive and equal to 
$Z\cdot\mathrm{SL}_{n}(K)$ where the torus $Z\cong\mathbb{G}_{\times}$ is the group
of scalar matrices.

\subsection{Chevalley groups}\label{chev}

\begin{definition} Let $G$ be an algebraic group defined over a field $k$.
Then $G$ is called \emph{$k$-split} if it has a maximal $k$-split 
torus.
\end{definition}

There is a unique simple, simply connected and $\mathbb Q$-split algebraic group of any given Lie type $\mathcal X$ and this is the called the \emph{Chevalley group of type 
$\mathcal X$}. 
There is a simple conceptial way to define the adjoint group 
$\overline{G} = G/Z(G)$, as described for example in \cite{carter}, Chapter 1: 
As we have seen $\overline{G}$ acts faithfully on the Lie algebra $L=L(G)$ of $G$ and so can be identified with a subgroup 
of $\mathrm{Aut}(L)$. In fact $\overline G$ is defined as the subgroup of $\mathrm{Aut}(L)$ generated by the elements
\[ \exp (\mathrm{ad} (x))= 1 + \mathrm{ad} (x) + \frac{\mathrm{ad}(x)^2}{2!} + \frac{\mathrm{ad}(x)^3}{3!} + \cdots \]  
where $x$ is an element of a root subgroup of $L$. Note that for such $x$ the linear transformation 
$\mathrm{ad}(x): L \rightarrow L$ is nilpotent and so the above series is finite. 

Moreover as described in 
\cite{carter} one can find a Lie subring $K$ of $L$ which is a finitely generated torsion 
free $\mathbb Z$-lattice of $L$ and such that $\exp (\mathrm{ad}(x))$ stabilizes $K$ for each $x$ as above.
Hence $\overline{G}$ is in fact defined over $\mathbb Z$ and one sees that 
the same is true for the universal cover $G$. Therefore its $R$-rational points $G_R$ are defined for any ring $R$. 
In particular $G_\mathbb{F}$ is defined for any finite field $\mathbb F$. As we shall see in section 
\ref{simpleslie} this is the construction of the untwisted finite simple groups of Lie type.

\section{Arithmetic groups and the congruence topology}

In this section and below $k$ will denote a number field (a finite extension
of $\mathbb{Q}$) and $\mathcal{O}$ its ring of integers. By convention,
\emph{prime ideals} of $\mathcal{O}$ are assumed nonzero. We begin by
recalling some information about the ring $\mathcal{O}$.

\subsection{Rings of algebraic integers in number fields\label{algnum1}}

$\mathcal{O}$ is the collection of all elements $x\in k$ satisfying a
polynomial equation
\[
x^{n}+a_{1}x^{n-1}+\cdots+a_{n-1}x+a_{n}=0
\]
with leading coefficient 1 and each $a_{i}\in\mathbb{Z}$. This is in fact a
subring of $k$. As an additive group it is isomorphic to $\mathbb{Z}^{d}$, the
free abelian group of rank $d$, where $d=(k:\mathbb{Q})$.

The ring $\mathcal{O}$ has Krull dimension $1$, i.e. every prime ideal is
maximal. Moreover, every nonzero ideal has finite index in $\mathcal{O}$. Each
nonzero ideal $I$ can be factorized
\[
I=\mathfrak{p}_{1}^{e_{1}}\cdot\mathfrak{p}_{2}^{e_{2}}\cdot\cdots
\cdot\mathfrak{p}_{m}^{e_{m}}%
\]
as a product of prime ideals $\mathfrak{p}_{i}$ and this factorization is
unique up to reordering of the factors. The Chinese Remainder Theorem says
that then%
\[
\mathcal{O}/I\cong\mathcal{O}/\mathfrak{p}_{1}^{e_{1}}\oplus\mathcal{O}%
/\mathfrak{p}_{2}^{e_{2}}\oplus\cdots\oplus\mathcal{O}/\mathfrak{p}_{m}%
^{e_{m}}.
\]

Each prime ideal $\mathfrak{p}$ divides (i.e. contains) a unique rational
prime $p\in\mathbb{N}$, and then $\mathfrak{p}\cap\mathbb{Z}=p\mathbb{Z}$. The
quotient $\mathcal{O}/\mathfrak{p}$ is a finite field of characteristic $p$.

If $p\mathcal{O}=\prod_{i=1}^{g}\mathfrak{p}_{i}^{e_{i}}$ is the factorization
of the principal ideal $(p)=p\mathcal{O}$ then
\begin{equation}
d=(k:\mathbb{Q})=\sum_{i=1}^{g}e_{i}n_{i},\ \text{ where }\left\vert
O/\mathfrak{p}_{i}\right\vert =p^{n_{i}}. \label{sumformula}%
\end{equation}

If $k$ is a Galois extension of $\mathbb{Q}$ then $e_{1}=\cdots=e_{g}$ and
$n_{1}=\cdots=n_{g}$. Also $e_{i}\not =1$ for at most finitely many rational
primes $p$ (the so-called \emph{ramified} ones). The \textbf{Chebotarev density
theorem} (see \cite{PR}, Chapter 1) implies that for a positive proportion of all
rational primes $p$ we have $g= (k: \mathbb Q)$ and $n_1= \cdots n_g=1$, i.e. 
the ideal $(p)$ splits in $k$. More precisely 
\[ \lim_{ x \rightarrow \infty} \frac{ |\{ p \leq x \ | \ p \textrm{ a prime%
 which splits in } k \}| }{ |\{ p \leq x \ | \ p \textrm{ prime }\}|} = 
\frac{1}{ |\mathrm{Gal}(k/\mathbb Q)|}=\frac{1}{d}.\]

Here $\mathrm{Gal}(k/ \mathbb Q)$ is the Galois group of $k$ over $\mathbb Q$.
\medskip

Let $S$ be a finite set of prime ideals. An element $a\in k$ is said to be
$S$\emph{-integral }if $Ja\subseteq\mathcal{O}$ where $J$ is some product of
prime ideals in $S$. The set of all $S$-integral elements forms a subring
$\mathcal{O}_{S}$ of $k$, containing $\mathcal{O}$, called the ring of
$S$\emph{-integers} of $k$. Of course, when $S$
is empty $\mathcal{O}_{S}=\mathcal{O}$.

\subsection{The congruence topology on $\mathrm{GL}_{n}(k)$ and
$\mathrm{GL}_{n}(\mathcal{O})$}

The \emph{congruence topology} on $k$ has as base of open neighbourhoods of
$0$ the set of all nonzero ideals of $\mathcal{O}$. The congruence topology on
$M_{n}(k)=k^{n^{2}}$ is then the product topology, and the congruence topology
on $\mathrm{GL}_{n}(k)$ (and on any closed subgroup) is the one induced by
that on $M_{n}(k)$. This means that a base of neighbourhoods of $1$ is the set
of subgroups $\mathrm{GL}_{d}(k)\cap(1_{n}+M_{n}(I))$ with $I$ a nonzero ideal
of $\mathcal{O}$.

More generally, for any set $S$ of prime ideals, we define the $S$%
\emph{-congruence topology} by taking only ideals that are products of prime
ideals not in $S$; equivalently, we can take as neighbourhood basis the set of
all nonzero ideals of $\mathcal{O}_{S}$.

It is easy to see that the congruence topology on $k$ and hence on $M_n(k)$ is Hausdorff: if $x \not = y$ 
are two elements of $k$ then there is an ideal $I$ of $\mathcal O$ such that $ x-y \not \in I$ and hence
 $(x+I)  \cap (y+I) = \emptyset$. In fact the congruence topology on $M_{n}(k)$ is finer than the Zariski topology as the 
following proposition demonstrates.

\begin{proposition} Let $W$ be a $k$-defined Zariski closed set of $M_n(K)=K^{n^2}$ defined by an ideal 
$T$ of polynomials in its $n^2$ coordinates. Then $W_k$ is closed in the congruence topology of $M_n(k)$.
\end{proposition}

\textbf{Proof.} Let $\mathbf x \in M_n(k)$ be an element of the congruence closure of $W_k$. So for any ideal
 $I$ of $\mathcal O$ 
we have an element $\mathbf y \in W_k$ such that $\mathbf x \equiv \mathbf y$ mod $I$. 
Now let $p$ be a polynomial from $T$ with coefficients in $k$. 
We may assume that up to a scalar multiple $p$ has coefficients from $\mathcal O$. But then 
$p(\mathbf x)\equiv p(\mathbf y)=0$ mod $I$, so $p(\mathbf x) \in I$ for any ideal $I$ of $\mathcal O$. 
This is possible only if $p(\mathbf x)=0$. Since this holds for all polynomials $p$ with coefficients in $k$, 
and since $W$ is defined over $k$ we deduce that
$\mathbf x \in W$. \bigskip

\textbf{Example:} Let $k=\mathbb Q$. Then any finite union or intersection of sets of the form 
\[ \{a+m\mathbb Z\} \times \{ b+n \mathbb Z\} \subset \mathbb Q^2 \quad a,b,m,n \in \mathbb Z \]
is an open set in the congruence topology of $\mathbb Q^2$ but none of them is Zariski open.

\bigskip

It is thus clear that the congruence topology on $M_n(k)$ has many more closed sets than the Zariski topology. So a Zariski dense subset of 
matrices may be rather 'sparse' in the congruence topology. From this point of view it is indeed surprising that in the case of 
simple algebraic groups the property of being Zariski dense has rather strong implication for the congruence closure of a subgroup. 
This is the main content of Theorem \ref{nori} below.  

\subsubsection{Valuations of $k$}

For any prime ideal $\mathfrak{p}$ of $\mathcal{O}$ the $\mathfrak{p}$-adic
topology is defined in the same way as the congruence topology except that the
ideals are only allowed to be positive powers of $\mathfrak{p}$. The
completion of $k$ with respect to this topology is denoted $k_{\mathfrak{p}}$
and the closure of $\mathcal{O}$ in $k_{\mathfrak{p}}$ is denoted
$\mathcal{O}_{\mathfrak{p}}$.

The valuation $v_{\mathfrak{p}}$ on $k_{\mathfrak{p}}$ is defined by
$v_{\mathfrak{p}}(a)=t$ where $t\in\mathbb{Z}$ is the largest integer such
that $\mathfrak{p}^{-t}a\subseteq\mathcal{O}_{\mathfrak{p}}$ (if $a\neq0$; one
sets $v_{\mathfrak{p}}(0)=\infty$). Thus $\mathcal{O}_{\mathfrak{p}}$ is the
\emph{valuation ring}, consisting of all elements of $k$ having valuation
$\geq0$; this implies that $\mathcal{O}_{\mathfrak{p}}$ is a local ring,
having $\mathfrak{p}\mathcal{O}_{\mathfrak{p}}$ as its unique maximal ideal.
(One often associates to such a valuation $v_{\mathfrak{p}}$ the corresponding
\emph{absolute value}: $\left\vert a\right\vert _{\mathfrak{p}}%
=q^{-v_{\mathfrak{p}}(a)}$ where $q=\left\vert \mathcal{O}/\mathfrak{p}%
\right\vert $, which is multiplicative.)

\begin{example}
[The $p$-adic numbers]\label{qp} Take $k=\mathbb{Q}$ with ring of integers $\mathbb{Z}$.
Let $p$ be a prime. The $p$\emph{-adic} valuation $v_{p}(x)$ is the usual one
where $v_{p}(x)=t$ is the largest integer such that $x=p^{t}a/b$ with integers
$a$ and $b$ coprime to $p$. A base for neighbourhoods of $0$ in the $p$-adic
topology on $\mathbb{Q}$ is the family of subgroups $\{\frac{p^{l}a}%
{b}\ |\ a,b\in\mathbb{Z},(p,b)=1\}$, $l=1,2,\ldots$. The completion of
$\mathbb{Q}$ with respect to this topology is the field $\mathbb{Q}_{p}$ of
$p$\emph{-adic numbers}. Inside $\mathbb{Q}_{p}$ we have the closure
$\mathbb{Z}_{p}$ of $\mathbb{Z}$, which is the ring of $p$\emph{-adic
integers}. We can view $\mathbb{Z}_{p}$ as a ring of infinite power series in
$p$:
\[
a_{0}+a_{1}p+\cdots+a_{k}p^{k}+\cdots,\quad a_{i}\in\{0,1,\ldots,p-1\}
\]
with the obvious addition and multiplication. The finite such sums comprise
the subring $\mathbb{Z}$. The unique maximal ideal is $p\mathbb{Z}_{p}$ and
every element $x\in\mathbb{Q}_{p}$ can be written uniquely as $x=p^{t}y$ for
some $y\in\mathbb{Z}_{p}\smallsetminus p\mathbb{Z}_{p}$ and $t\in\mathbb{Z}$.
\end{example}

In general if ${\mathfrak{p}}_{i}$ is a prime ideal of $\mathcal{O}$ dividing
$p$ then $k_{{\mathfrak{p}}_{i}}$ is a vector space over $\mathbb{Q}_{p}$ of
dimension $e_{i}n_{i},$ the number in (\ref{sumformula}) above.

The valuations $v_{\mathfrak{p}}$ for all prime ideals $\mathfrak{p}$ of
$\mathcal{O}$ comprise the \emph{non-archimedean} valuations of $k$. Now,
suppose that the number field $k$ has $s$ embeddings $v_{i}:k\rightarrow
\mathbb{R}$ ($i=1,\ldots,s$) and $2t$ non-real embeddings $v_{j},~\overline
{v}_{j}:k\rightarrow\mathbb{C}$ ($j=s+1,\ldots,s+t$). Composing these with the
ordinary real or complex absolute value gives the set $V_{\infty}$ of $s+t$
\emph{archimedean} absolute values on $k$. For $v\in V_{\infty}$ we put
$k_{v}=\mathbb{R}$ or $k_{v}=\mathbb{C}$ according as the corresponding
embedding of $k$ is real or non-real.

The ring of $S$-integers has a more natural definition in terms of
valuations:
\[
\mathcal{O}_{S}=k\cap%
{\displaystyle\bigcap\limits_{\mathfrak{p}\notin S}}
\mathcal{O}_{\mathfrak{p}}.
\]

A word of warning: the notation $\mathcal O_S$ can be a bit confusing: if $S=\{\mathfrak{q}\}$ consists of a single prime then $\mathcal{O}_{\{\mathfrak{q}\}}$
is the ring of $\{\mathfrak{q}\}$-integers, while $\mathcal{O}_{\mathfrak{q}}$ is
the completion of $\mathcal O$ at $\mathfrak{q}$.  For example $\mathbb{Z}_{\{p\}}=
\mathbb{Z}[1/p]\subset \mathbb{Q}$ while $\mathbb{Z}_p$ is the ring of $p$-adic
integers.

\subsection{Arithmetic groups}

Suppose we are given a linear algebraic group $G$ defined over $k$ with a faithful
representation $G\hookrightarrow\mathrm{GL}_{n}(K)$, also defined over $k$.

\begin{definition}
A subgroup $\Gamma$ of $G_{k}$ is called \emph{arithmetic} if it is
commensurable with the group of $\mathcal{O}$-integral points $G_{\mathcal{O}%
}$ (in other words $\Gamma\cap G_{\mathcal{O}}$ has finite index in both
$\Gamma$ and $G_{\mathcal{O}}$).
\end{definition}

It turns out that this definition is independent of the choice of $k$-defined
linear representation of $G$.

More generally we can define the $S$\emph{-arithmetic subgroups} of $G(k)$ as
those commensurable with $G_{\mathcal{O}_{S}}$. When the set $S$ has not been
specified we shall always assume that it is empty.

The simplest examples of arithmetic groups are $(\mathcal{O},+)$ and
$(\mathcal{O}^{\ast},\times)$ the additive and multiplicative groups of the
ring of integers of $k$. We thus see that the study of arithmetic groups is a
generalization of classical algebraic number theory.

One of the most general results about arithmetic groups is the following

\begin{theorem}
\emph{(\cite{PR}, chapter 4)}\label{FP}  Let $\Gamma$ be an arithmetic subgroup of a $k$-defined
linear algebraic group $G$ as above. Then $\Gamma$ is finitely presented and
has only finitely many conjugacy classes of finite subgroups.
\end{theorem}

For $S$-arithmetic groups the above statement is still true, provided that $G$
is reductive. \medskip

Now an ($S$-)arithmetic group $\Gamma$ has its own ($S$-)congruence topology
induced from the ($S$-)congruence topology of $\mathrm{GL}_{n}(k)$. We call a
subgroup $\Delta\leq\Gamma$ an \ ($S$-)congruence subgroup if is is open in this
topology, i.e. if $\Delta$ contains a \emph{principal congruence subgroup}
$\Gamma\cap(1_{n}+M_{n}(I))$ for some nonzero ideal $I$ of (coprime to $S$).
The \emph{congruence images} $\Gamma/N$ of $\Gamma$ are those with kernel a
congruence subgroup $N\vartriangleleft\Gamma$.

Clearly a congruence subgroup of $\Gamma$ has finite index, but the converse
is not true in general. When it does hold, that is if every subgroup of finite
index is a congruence subgroup, $\Gamma$ is said to have \emph{the congruence
subgroup property} (CSP).

There is a neat way to state CSP in term of profinite groups. If $\mathcal{X}
$ is an intersection-closed family of normal subgroups of finite index in
$\Gamma$, one defines the $\mathcal{X}$-\emph{completion} of $\Gamma$ to be
the inverse limit%
\begin{align*}
\widehat{\Gamma}_{\mathcal{X}}  &  =\underset{\longleftarrow}{\lim}%
_{N\in\mathcal{X}}\Gamma/N\\
&  =\left\{  (\gamma_{N})_{N\in\mathcal{X}}\mid p_{NM}(\gamma_{N})=\gamma
_{M}~\forall N\leq M\in\mathcal{X}\right\}  \leq%
{\displaystyle\prod_{N\in\mathcal{X}}}
\Gamma/N,
\end{align*}
where $p_{NM}:\Gamma/N\rightarrow\Gamma/M$ denotes the natural quotient map
for each $N\leq M$. (With the topology induced from the product topology on
the Cartesian product, $\widehat{\Gamma}_{\mathcal{X}}$ becomes a compact
topological group, a \emph{profinite group}). 

A natural example of inverse limits is the valuation ring $\mathcal{O}_{\mathfrak{p}}$: For a prime ideal $\mathfrak{p}$ of $\mathcal O$ the inverse
limit
\[ \underset{\longleftarrow}{\lim}_{n \in \mathbb{N}}  \mathcal O/\mathfrak{p}^n \mathcal O \]
is isomorphic as a ring to the completion $\mathcal{O}_{\mathfrak{p}}$ of
$\mathcal O$ with respect to the $p$-adic topology defined by the powers of the ideal $\mathfrak{p}$. This also shows that $\mathcal{O}_{\mathfrak{p}} / \mathfrak{p}^n \mathcal{O}_{\mathfrak{p}}$ is isomorphic to $\mathcal O/ \mathfrak{p}^n
\mathcal{O}.$

\medskip

We are interested in two special
choices for $\mathcal{X}$. When $\mathcal{X}$ consists of \emph{all} normal
subgroups of finite index, $\widehat{\Gamma}_{\mathcal{X}}=\widehat{\Gamma}$
is the \emph{profinite completion} of $\Gamma$. When $\mathcal{X}$ consists of
all the normal \emph{congruence} subgroups, $\widehat{\Gamma}_{\mathcal{X}%
}=\widetilde{\Gamma}$ is the \emph{congruence completion} of $\Gamma$. There
is an obvious natural projection $\pi:\widehat{\Gamma}\rightarrow
\widetilde{\Gamma}$, which is clearly surjective.

Now we can reformulate the congruence subgroup property as saying that the map
$\pi$ is bijective. For many purposes the following generalization of CSP is
more relevant: the arithmetic group $\Gamma$ is said to have the
\textit{generalized congruence subgroup property} (GCSP for short) if the
kernel of $\pi:\widehat{\Gamma}\rightarrow\widetilde{\Gamma}$ is finite. Group
theoretically this says that any subgroup of finite index in $\Gamma$ is
commensurable `with bounded index' with a congruence subgroup. There is a
famous conjecture by Serre which characterizes the $S$-arithmetic groups (in
semisimple algebraic groups) with GCSP as those having $S$-rank at least 2:
see section \ref{serre}.

\section{The Strong Approximation Theorem}

The congruence images of the $S$-arithmetic group $\Gamma=G_{\mathcal{O}_{S}}$
are easier to understand when $G$ has the \textit{strong approximation
property}. In order to explain this we need several more definitions. \medskip

Recall that $k_{\mathfrak{p}}$ and $\mathcal{O}_{\mathfrak{p}}$ are the
completions of $k$ and $\mathcal{O}$ with respect to the $\mathfrak{p}$-adic
topology defined by powers of the prime ideal $\mathfrak{p}\vartriangleleft
\mathcal{O}$. As usual we set $G_{k_{\mathfrak{p}}}=G\cap M_{n}%
(k_{\mathfrak{p}})$ and $G_{\mathcal{O}_{\mathfrak{p}}}=G\cap M_{n}%
(\mathcal{O}_{\mathfrak{p}})$. The first of these is a locally compact totally
disconnected topological group and the second is a compact subgroup. In fact
$G_{\mathcal{O}_{\mathfrak{p}}}$ is an example of a $p$\textit{-adic analytic group}. We
refer to $G_{k_{\mathfrak{p}}}$ as the completion of $G$ at $\mathfrak{p}$.

Similarly, if $v$ is an archimedean real (resp. complex) absolute value of $k
$ associated to an embedding $\nu_{i}$, then we write $G_{v}$ for
$G_{\mathbb{R}}^{\nu_{i}}$ (resp. $G_{\mathbb{C}}^{\nu_{i}}$) where
$G^{\nu}\leq\mathrm{GL}_{n}(\mathbb{C})$ is the group obtained by applying
$\nu$ to the defining equations of the affine variety $G$.

The profinite groups $G_{\mathcal O_{\mathfrak p}}$ are in close relationship
with the congruence
images of $G_{\mathcal O_S}$:

Recall that the algebraic group $G$ is $k$-defined. It is easy to see that for almost all prime ideals $\mathfrak{p}$ the coefficients of the equations
defining $G$ in $\mathrm{GL}_n$ are not divisible by $\mathfrak p$. Therefore
we can consider these equations modulo $\mathfrak p^n$ for any $n \in \mathbb{N}$. Denote the set of their solutions in
$\mathcal O/\mathfrak{p}^n$ by $G_{\mathcal O/\mathfrak{p}^n}$:
this is a finite subgroup of $\mathrm{GL}_n(\mathcal O/\mathfrak{p}^n)$ and
is called the \emph{reduction} of $G$ modulo $\mathfrak p^n$. See \cite{PR} p. 142-146 for more details about reductions of affine algebraic varieties
and groups.

Now consider the quotient mapping
\[ \mathcal{O}_{\mathfrak p} \rightarrow \mathcal{O}_{\mathfrak{p}} / 
\mathfrak{p}^n \mathcal{O}_{\mathfrak{p}} \simeq \mathcal O/ \mathfrak{p}^n 
\mathcal{O}.\] 

This induces a homomorphism 
\[ \pi_{\mathfrak{p}^n}: \  G_{\mathcal O_{\mathfrak{p}}} \rightarrow G_{\mathcal{O}/ \mathfrak{p}^n \mathcal O}.\] 

The following is the content of Proposition 3.20 of \cite{PR}. 
\begin{proposition} \label{redp} The
maps $\pi_{\mathfrak{p}^n}$ are surjective for all but finitely many primes
$\mathfrak{p}$ (and all integers $n$).
\end{proposition}
We say that $G$ has \emph{good reduction} for such primes $p$. \medskip

Assume from now on that $\mathfrak{p}$ is not in the finite set $S$. The restriction of $\pi_{\mathfrak{p}^n}$ to its dense subgroup $G_{\mathcal O_S} \leq  
G_{\mathcal O_{\mathfrak{p}}}$ is the homomorphism

\[ G_{\mathcal O_{S}} \rightarrow G_{\mathcal{O}_S/ \mathfrak{p}^n \mathcal O_S} \] obtained by reducing all entries of $\Gamma=G_{\mathcal{O}_S} \leq
\mathrm{GL}_n(\mathcal O_S)$ modulo $\mathfrak p^n$. So the images of $\pi_{\mathfrak{p}^n}$ are all congruence images of
$\Gamma$. What is not clear at this point is how to combine these
to describe the congruence
images of $\Gamma$ at composite ideals. This is the content of the strong
approximation theorem below.

\bigskip

Define
\[
G_{S}:=\prod_{v\in V_{\infty}}G_{v}\times\prod_{p\in S}G_{k_{p}}.
\]
This is a locally compact group and the image of $\Gamma$ in $G_{S}$ under the
diagonal embedding in each factor is a \textit{lattice} in $G_{S}$, i.e., a
discrete subgroup of finite co-volume. As a consequence the arithmetic
subgroup $\Gamma=G_{\mathcal{O}_{S}}$ is infinite if and only if the group
$G_{S}$ is non-compact.
\medskip

Let
\[
G_{\widehat{\mathcal{O}}_{S}} = \prod_{p \not \in S} G_{\mathcal{O}_{p}}.
\]

Again there is an obvious diagonal embedding $i:\Gamma\rightarrow
G_{\widehat{\mathcal{O}}_{S}}$ and the congruence topology of $\Gamma$
coincides with the topology induced in $i(\Gamma)$ as a subgroup of the
profinite group $G_{\widehat{\mathcal{O}}_{S}}$. Hence the congruence
completion $\widetilde{\Gamma}$ is isomorphic to the closure $\overline{i(G)}$
of $i(G)$ in $G_{\widehat{\mathcal{O}}_{S}}$. The strong approximation theorem
states that under certain conditions $i(G)$ is \emph{dense} in $G_{\widehat{\mathcal{O}}_{S}}$, 
and therefore
$\widetilde{\Gamma}\simeq G_{\widehat{\mathcal{O}}_{S}}$.

\begin{theorem}
[Strong approximation for arithmetic groups]\label{strongapprox}
\emph{(\cite{PR}, Theorem 7.12). \ }Let $G$ be a connected simple simply
connected algebraic group defined over a number field $k$ and let the groups
$\Gamma=G_{\mathcal{O}_{S}}$, $G_{S},G_{\widehat{\mathcal{O}}_{S}}$ and the
embedding $i:\Gamma\rightarrow G_{\widehat{\mathcal{O}}_{S}}$ be as above.
Assume that $\Gamma$ is infinite (which is equivalent to $G_{S}$ being
non-compact). Then $i(\Gamma)$ is dense in $G_{\widehat{\mathcal{O}}_{S}}$ and
hence $\widetilde{\Gamma}\simeq G_{\widehat{\mathcal{O}}_{S}}$.
\end{theorem}

\noindent When the conclusion holds we say that $G_{\mathcal{O}_{S}}$ has the 
\textbf{strong approximation property}, or that $G$ has the strong approximation property
w.r.t. $S$.
\medskip

Note: Usually the strong approximation theorem is formulated for the group of
$k$-rational points $G_{k}$ and says that $G_{k}$ is dense in the adelic group
$G_{\mathcal{A}_{S}}$, the statement we have given above is equivalent to
this (and more transparent for arithmetic groups); see \cite{PR}, Chapter 7.

More generally, a connected algebraic group $G$ has the strong approximation
property if its maximal reductive quotient $H=G/R_{u}(G)$ is a direct product
of simple simply connected groups, and $H_{S}$ is non-compact. \medskip

The strong approximation theorem can be viewed as a generalization of the
Chinese remainder theorem, which in this setting says that the diagonally
embedded image of $\mathbb{Z}$ is dense in $\prod_{p \ \mathrm{prime}}
\mathbb{Z}_{p}$. In the general situation the theorem says that the finite images of the product 
$G_{\widehat{\mathcal{O}}_{S}}$ coincide with the congruence images of $\Gamma$
\medskip

Note: The condition that $G$ be simply connected is indeed necessary (Exercise 10).

\medskip

 Set $\mathbb{F}_{q(\mathfrak{p})}$ where
$q(\mathfrak{p})=\left\vert \mathcal{O}/\mathfrak{p}\right\vert $.
Then theorem \ref{strongapprox} and Proposition \ref{redp} give

\begin{corollary} \label{qsimpleim}
Under the hypotheses of Theorem \ref{strongapprox}, we have $\pi
_{\mathfrak{p}}(\Gamma)=G_{\mathbb{F}_{q(\mathfrak{p})}}$ for all but finitely
many primes $\mathfrak{p}\notin S$.
\end{corollary}

In turn the groups $G_{\mathbb{F}_{q(\mathfrak{p})}}$ are easy to describe
when $G$ is semisimple, see Proposition \ref{imagep} below.

For the moment we shall note case the relationship between $G$ and $H$ when $G=\mathcal R_{k/\mathbb Q}H$ is a restiction of scalars of $H$.
\begin{proposition} \label{restr} Let $k$ be a finite Galois extension of $\mathbb Q$ and $G=\mathcal R_{k/\mathbb Q}H$ be the restriction of scalars of some $\mathbb Q$-defined algebraic group $H$. Then for almost all primes $p$ 
\[ G_{\mathbb F_p} = \prod_{\mathfrak{p}|p} H_{\mathbb F_{q(\mathfrak{p})}},\]
where the product on the right is over all prime ideals $\mathfrak{p}$ of $k$ dividing $p$.

\end{proposition} 

Indeed all but finitely many primes $p$ are unramified in $k$ and therefore 
\[\mathcal O /p\mathcal O= \prod_{\mathfrak{p}|p}\mathcal O/\mathfrak{p}=\prod_{\mathfrak{p}|p}\mathbb F_{q(\mathfrak{p})},\]
where $\mathcal O$ is the ring of integers of $k$.
The proposition follows immediately from $G_\mathbb Z= H_{\mathcal O}$ and Proposition \ref{redp}.
\subsection{An Aside: Serre's Conjecture}

\label{serre} We now have most of the definitions to state Serre's conjecture.

\begin{definition}
For a valuation $v$ of $k$ the $k_{v}$-rank of the topological group $G_{v}$
is the largest integer $n$ such that $G_{v}$ contains the direct product
$(k_{v}^{\ast})^{n}$. The $S$-rank of an algebraic group $G$ is
\[
\sum_{v\in V_{\infty}\cup S}k_{v}\text{-rank of }G_{v}%
\]
where $V_{\infty}$ is the set of all archimedean valuations of $k$.
\end{definition}

\begin{conjecture}
\emph{(J-P. Serre)}  A connected simply connected simple algebraic group $G$ has
the generalized $S$-congruence subgroup property if and only if the $S$-rank
of $G$ is at least $2$.
\end{conjecture}

For example the group $\mathrm{SL}_{n}(\mathbb{Z})$ has CSP if $n>2$ but not
if $n=2$.

Currently Serre's conjecture is open for some groups of $S$-rank 1 and also
when $G$ is a totally anisotropic form of $A_{n}$, see \cite{PR}, \S 9.5.

\section{The Nori-Weisfeiler theorem and Lubotzky's alternative}

It will be too much to expect that the Strong Approximation Theorem holds for
linear groups in general, indeed it doesn't hold for algebraic tori.
Nevertheless there is something that can be said when the group is non-soluble.

\begin{theorem}
\emph{(Nori \cite{nor}, Weisfeiler \cite{Weis})}\label{nori}  Let $\Delta$ be a
Zariski-dense subgroup of a $\mathbb{Q}$-simple simply connected linear
algebraic group $G\leq \mathrm{GL}_{n}(\mathbb{C})$ and suppose that
$\Delta\leq G_{\mathbb{Z}_{S}}$ for some finite set of primes $S$. Let
$i:\Delta\rightarrow G_{\widehat{\mathbb{Z}}_{S}}$ be the diagonal embedding.

Then the closure $\overline{i(\Delta)}$ of $i(\Delta)$ in $G_{\widehat
{\mathbb{Z}}_{S}}$ is an open subgroup of $G_{\widehat{\mathbb{Z}}_{S}}$.
\end{theorem}

It follows that for all but finitely many primes $p$, all the groups $G_{\mathbb{Z
}/(p^{m} \mathbb{Z})}$ appear as congruence images of $\Delta$. \bigskip

There are several different proofs of this theorem. We shall sketch one of
them in section \ref{nor}. For the moment, let us assume this result and
 deduce Theorem \ref{lubalt}. We restate it here: \medskip

\textbf{Theorem 1} \emph{Let $\Delta\leq \mathrm{GL}_{n}(k)$ be a finitely
generated linear group over a field $k$ of characteristic 0. Then one of the
following holds:} \medskip

 (a) \emph{the group $\Delta$ is virtually soluble, or \medskip}

 (b) \emph{there exist a connected simply connected $\mathbb{Q}$-simple
algebraic group $G$, a finite set of primes $S$ such that 
$\Gamma =G_{\mathbb{Z}_{S}}$ is infinite, and a subgroup $\Delta_{1}$ of finite index
in $\Delta$ such that every congruence image of $\Gamma$ appears as a quotient
of $\Delta_{1}$}. \bigskip

\textbf{Proof of Theorem 1: }Suppose that we have a finitely generated linear
group $\Delta\leq\mathrm{GL}_{n}(\mathbb{C})$. Then in fact 
$\Delta\leq\mathrm{GL}_{n}(J)$ for some finitely generated subring $J$ of $\mathbb{C}$.

Now the Jacobson radical (the intersection of the maximal ideals of $J$) is
trivial and so $J$ is residually a number field. Indeed if $\mathbf m$ is a maximal
ideal of $J$ then $J/\mathbf m$ is a finitely generated algebra which is a field. By
Corollary 7.10 in \cite{AM} ('The weak Nullstellensatz'), $J/\mathbf m$ is a finite
extension of $\mathbb{Q}$, i.e. a number field.

Hence $\Delta$ is residually in $\mathrm{GL}_{n}(k_{i})$ for some number fields $k_{i}$.
Suppose that $\Delta$ is not virtually soluble. By Lemma \ref{plat} it follows
that there is $i\in I$ such that the image of $\Delta$ in 
$\mathrm{GL}_{n}(k_{i})$ is not virtually soluble. Replacing $\Delta$ with this image we
may assume that $\Delta\leq\mathrm{GL}_{n}(k)$ for some number field $k$.

Consider $\mathrm{GL}_{n}(k)$ as a subgroup of $\mathrm{GL}_{nd}(\mathbb{Q})$
where $d=(k:\mathbb{Q})$. Let $\mathcal{G}$ be the Zariski-closure of $\Delta$
in $\mathrm{GL}_{nd}(K)$. This is a $\mathbb{Q}$-defined linear algebraic
group and we take its connected component $\mathcal{G}_{0}$ at the identity.

Let $\Delta_{1}=\mathcal{G}_{0}\cap\Delta$. This has finite index in $\Delta$
and is Zariski-dense in $\mathcal{G}_{0}$. Since $\Delta$ is not virtually
soluble the connected algebraic group $\mathcal{G}_{0}$ is not soluble. By
Exercise 12 we see that there exists a $\mathbb{Q}$-simple connected algebraic
group $G$ and a $\mathbb{Q}$-defined epimorphism $f:\mathcal{G}_{0}\rightarrow
G$. Now $f(\Delta_{1})$ is dense in $G$ and we may replace $\Delta$ by
$f(\Delta_{1})$ and $\mathcal{G}_{0}$ by $G$ to reduce the situation to where
we have a finitely generated Zariski-dense subgroup $\Delta\leq G_{\mathbb{Q}%
}$ of a $\mathbb{Q}$-simple connected linear algebraic group $G$. The main
difference with the setup of Theorem \ref{nori} is that $G$ may not be simply
connected. However $G$ is isogenous to its simply connected cover
$\widetilde{G}$, i.e., there is a $\mathbb{Q}$-defined surjection
$\pi:\widetilde{G}\rightarrow G$, where $\ker\pi=Z$ is a finite central
subgroup of $\widetilde{G}$.

It is not in general true that $\pi (\widetilde G_\mathbb Q)= G_\mathbb{Q}$
but at least we have the following

\begin{proposition} The group $G_\mathbb{Q}/\pi (\widetilde G_\mathbb Q)$
is abelian of finite exponent dividing $|Z|$.
\end{proposition}

\textbf{Proof:} Let $A$ be the Galois group of $K/\mathbb Q$ where $K$ is the algebraic closure of $\mathbb Q$. 
Then $\widetilde G_{\mathbb Q}$
consists of all those $g \in \widetilde G_K$ such that $g^\alpha=g$ for all $\alpha
\in A$. On the other hand $\pi^{-1} (G_\mathbb Q)$ consists of those $g \in \widetilde{
G}_K$ such that $g^\alpha \equiv g$ mod $Z$ for all $\alpha \in A$. Suppose that $g, h \in \pi^{-1} (G_\mathbb Q)$, thus $g^\alpha \equiv g $ and $h^\alpha
\equiv h$ mod $Z$ for all $\alpha \in A$. Now using
that $Z$ is central in $\widetilde G$ we see that 
$[g,h]^\alpha= [g^\alpha,h^\alpha]=[g,h]$ and hence that $[g,h] \in
\widetilde G_{\mathbb Q}$. Let $m= \exp Z$. In the same way also we see that if $g^\alpha
 \equiv g$ mod $Z$ then $(g^m)^\alpha= g^m$ and therefore $g^m \in \widetilde G_{\mathbb Q}$. So $\pi^{-1} (G_\mathbb Q)/\widetilde G_{\mathbb Q}$ is abelian
of exponent dividing $|Z|$ and this implies the Proposition. $\square$
\medskip

Now take $\Delta_{0}=\Delta\cap\pi(\widetilde G_{\mathbb{Q}})$, this is a
subgroup of finite index in $\Delta$ because $\Delta/ \Delta_{0}$ is a
finitely generated abelian group of finite exponent. Let $U_{0} = \pi
^{-1}(\Delta_{0}) \cap \widetilde G_{\mathbb{Q}}$, then $U_0/(U_0 \cap Z) \simeq \Delta_0$; $U_{0}$ is a
finitely generated linear group it is residually finite. So we can find a
subgroup $U$ of finite index in $U_{0}$ such that $U \cap Z=\{1\}$.
Then
$U$ is isomorphic to $\pi(U)$ which is a subgroup of finite index in
$\Delta_{0}$ and hence in $\Delta$.

Now take $\Delta_{1}= \pi(U) \simeq U$. Observe that $U$ is Zariski dense in
the $\mathbb{Q}$-simple, connected and simply connected algebraic group
$\widetilde G$. In addition $U$ is finitely generated and inside $\widetilde
G_{\mathbb{Q}}$. It follows that there is a finite set $S$ of
rational primes such that $U \leq\widetilde G_{\mathbb{Z}_{S}}$. All the
conditions of Theorem \ref{nori} are now satisfied with $U$ and $\widetilde G$ in place of $\Delta$ and $G$.
Hence we deduce that the congruence completion of $U$ is an open subgroup of
\[
G_{S} = \prod_{p \not \in S} G_{\mathbb{Z}_{p}}.
\]
This open subgroup projects onto all but finitely many of the factors in the product $G_{S}$. So
by enlarging $S$ to some finite set $S_1$ we may ensure that the congruence
completion of $U$ maps onto $\prod_{ p \not \in S_{1} }G_{\mathbb{Z}_{p}}$.
Since $U$ is isomorphic to $\Delta_{1}$ Theorem \ref{lubalt} follows.

\section{Some applications to Lubotzky's alternative}

\label{applic}

As noted in the introduction, Theorem \ref{lubalt} puts a substantial
restriction on the finite images of a linear group in characteristic $0$.
First we need to introduce

\subsection{The finite simple groups of Lie type.}\label{simpleslie}

For a detailed account of the material of this section we refer to Carter's
book \cite{carter}.

The \emph{untwisted} simple groups of Lie type are the groups $L=G_{\mathbb{F}_{q}}/Z$
where $G$ is a simply connected Chevalley group defined over $\mathbb Z$ and $Z$ is the centre of the group of 
rational points $G_{\mathbb{F}_{q}}$ over the finite field $\mathbb{F}_{q}$. The \emph{type} of $L$ is
just the Lie type $\mathcal{X}$ of $G$.

The \emph{twisted} simple groups arise as the fixed points $L^{\sigma}$ of a
specific automorphism $\sigma$ (of order 2 or 3) of some untwisted simple group $L$.
Such twisted Lie type simple groups are for example $\mathrm{PSU}_{n}(q)$. The
(untwisted) type of $L^{\sigma}$ is just the Lie type of $L$. For example the
untwisted Lie type of $\mathrm{PSU}_{n}(q)$ is $A_{n-1}$. \bigskip

A finite group $L$ is quasisimple if $L=[L,L]$ and $L/Z(L)$ is simple. Similarly
to the isogenies described in Theorem \ref{K-simple}, the quasisimple finite groups break up into
families with the same simple quotient. The members of each family have the same simple quotient, say $S$ and there is a largest member of the family
$L$, called \emph{the universal cover} of $S$. All the other members of
the family are the quotients $L/A$ where $A \leq Z(S)$. The type (twisted
or not) of a quasisimple group is the same as that of its simple quotient.

\subsection{Refinements}

Let us return to Corollary \ref{qsimpleim}. Recall that $G$ was a simple
simply connected linear algebraic group defined over an algebraic number
field $k$ and $\Gamma = G_{\mathcal O_S}$ for a ring of algebraic $S$-integers $\mathcal O_S$ of $k$
. The group $\Gamma$ then maps onto $G_{ \mathbb{F}_{q(\mathfrak{p})}}$ for almost all
$\mathfrak{p} \not \in S$. 

\begin{proposition} \label{imagep} Assume in the above situation that $G$ is
absolutely simple. Then for almost all prime ideals $\mathfrak{p}$ outside $S$ the reduction $G_{\mathbb{F}_{q(\mathfrak{p})}}$ of $G$ modulo $\mathfrak{p}$
is a 
quasisimple finite group. 
\end{proposition}

Now from the description of the
$k$-forms of $G$ it follows that $G$ splits over $\mathbb F_{q(\mathfrak p)}$ if and only if some 
specific polynomials in $k[x]$ (depending only on $G$) splits completely in linear factors in the finite field 
$\mathbb{F}_{q(\mathfrak p)}$. The \textbf{Chebotarev density Theorem} now implies that 
$G_{\mathbb{F}_{q(\mathfrak p)}}$  is an \textbf{untwisted} quasisimple group
for a positive proportion of the primes $\mathfrak p$ of $k$.
\medskip

Next consider the situation of Theorem \ref{lubalt}. There we have a $\mathbb Q$-simple algebraic group $G$ such that all congruence images of $G_{\mathbb{Z}_S}$ occur as quotients of $\Delta_1$. Now $G$ may not be absolutely simple but in any case there is a finite Galois extension $k$ of $\mathbb Q$ and an absolutely simple $k$-defined group $H$ such that $G=\mathcal R_{k/\mathbb Q} H$. Proposition \ref{restr} gives that for almost all rational primes $p$ outside $S$
\[ G_{\mathbb F_p} = \prod_{\mathfrak{p}|p} H_{\mathbb F_{q(\mathfrak{p})}}.\]
where as before the product on the right is over all prime ideals $\mathfrak{p}$ of $k$ dividing the (unramified) prime $p$. Note that the degree of $\mathbb F_{q(\mathfrak p)}$ over $\mathbb F_p$ is bounded by $(k:\mathbb Q)$

Therefore Theorem \ref{lubalt} in combination with Corollary \ref{qsimpleim} gives

\begin{corollary}
Suppose that $\Gamma< \mathrm{GL}_{n}(K)$ is a finitely generated linear
group in characteristic 0 which is not virtually soluble. Then there is

\begin{itemize}
\item a positive integer $d$,

\item a Lie type $\mathcal{X}$,

\item for every prime $p$ a finite field $\mathbb F_{p^f}$ of degree $f \leq d$ over $\mathbb F_p$ and a finite simple group $L(p^f)$ of Lie type over
$\mathbb{F}_{p^f}$ whose \emph{untwisted} type is $\mathcal{X}$ (e.g. if
$\mathcal{X}=A_{n-1}$ then $L(p^f)$ is either $\mathrm{PSL}_{n}(p^f)$ or $\mathrm{PSU}_{n}(p^f)$), and

\item a subgroup of finite index $\Gamma_{0}$ in $\Gamma$,
\end{itemize}

such that $\Gamma_{0}$ maps onto $L(p^f)$ for almost all primes $p$. Moreover,
for a positive proportion of these primes one has $f=1$ and the group $L(p)$ is untwisted.
\end{corollary}

One consequence of this is that $\Gamma$ cannot have polynomial subgroup
growth because the Cartesian product $\prod_{p \ \mathrm{prime}}L(p)$ doesn't
have polynomial subgroup growth, see \cite{SG} Chapter 5.2 for details.
\bigskip

The untwisted type $\mathcal{X}$ of the simple groups $L(p)$ is not
completely arbitrary: Let $G$ be the simple algebraic group of type
$\mathcal{X}$ as stated in Theorem \ref{K-simple}. Then $G$ is an image of the
connected component of the Zariski closure of $\Gamma$ in $G_{n}(K)$.

There is one particular case when the group $G$ is explicitly determined: when
$\Gamma$ is a subgroup of $\mathrm{GL}_{2}(\mathbb{C})$. Then the dimension of
$G$ is at most 4. On the other hand from the classification in Theorem
\ref{K-simple} it follows that the only simple algebraic group of dimension
less than 8 is $SL_{2}$. Therefore we obtain the following

\begin{proposition}
A finitely generated subgroup $\Gamma$ of $\mathrm{GL}_{2}(\mathbb{C})$ which
is not virtually soluble has a subgroup of finite index $\Gamma_{0}$ which
maps onto $\mathrm{PSL}_{2}(p)$ for infinitely many, in fact for a positive proportion
of all primes $p$.
\end{proposition}

This result is used in \cite{LLR} where the authors prove that any lattice
$\Lambda$ in $\mathrm{PSL}_{2}(\mathbb{C})$ has a collection $\{N_{i}\}_{i}$ of
subgroups of finite index such that $\bigcap_{i} N_{i}=\{1\}$ and $\Lambda$
has property $\tau$ with respect to $\{N_{i}\}_{i}$. As a corollary the
authors obtain that any hyperbolic 3-manifold has a co-final sequence of
finite covers with positive infimal Heegaard gradient.

\subsection{Normal subgroups of linear groups}

A normal subgroup $N$ of a finitely generated group of course does not need to be finitely generated. So it comes as no surprise that when this happens
in linear groups we can put further restriction of the finite images of $N$.

\begin{proposition} \label{norm} Let $\Gamma$ be a finitely generated linear group with
a finitely generated normal subgroup $\Delta$. Assume that $\Delta$ is not
virtually soluble. Then there exist a number $C>0$ and a Lie type
$\mathcal X$ with the following property: For infinitely
many primes $p \in \mathbb N$ the group $\Delta$ has a normal $\Gamma$-invariant
subgroup $N$ with $\Delta/N$ isomorphic to a direct product of at most $k$
copies of the untwisted finite simple group $L(p)$ of type $\mathcal X$ over $\mathbb
F_p$. 
\end{proposition}

\textbf{Sketch of Proof:} Using similar arguments to those in the proof of Theorem \ref{lubalt} we can reduce to the case 
when $\Gamma \leq \mathrm{GL}_d(\mathbb Q)$ for some integer $d$ and $\Delta$ is Zariski dense in  some absolutely semisimple simply connected algebraic group $G \leq  \mathrm{GL}_d$
defined over $\mathbb Q$ with isomorphic simple factors. Moreover we have $\Gamma \leq \mathrm{GL}_d(\mathbb Z_S)$ for some finite set of rational primes $S$.

Let $t$ be the number of simple factors of $G$.

As before for a rational prime $p \not \in S$ let $\pi_p$ be the homomorphism 
$\mathrm{GL}_d(\mathbb Z_S) \rightarrow \mathrm{GL}_d(\mathbb
F_p)$ obtained by reducing $\mathbb Z_S$ mod $p$.

From Theorem \ref{nori} we deduce that for all but finitely many primes $p$
outside $S$ one has $\pi_p(\Delta)=G_{\mathbb F_p}$= $\pi_p (G_{\mathbb
Z_S})$. Let $M_p= \ker \pi_p$ and
$N_p= \Delta \cap M_p$. Then $\Delta /N_p \simeq G_{\mathbb F_p}$ is a central
product of at most $t$ quasisimple groups of the same Lie type as the factors
of $G$. Also for infinitely many primes $p$ these factors are untwisted quasisimple
groups.

Now the only thing remaining is to observe that $N_p$ is normal in
$\mathrm{GL}_d(\mathbb Z_S)$ and therefore $N_p=\Delta \cap M_p$ is invariant under
$\Gamma$. Hence $G_{\mathbb F_p}/ \mathrm{Z}(G_{\mathbb F_p})$ is the required $\Gamma$-invariant
quotient of $\Delta$.
\bigskip

As suggested b Lubotzky Proposition \ref{norm} may be relevant in the following open problem: 

\begin{conjecture} Let $n>2$ and consider $\mathrm{Aut}(F_n)$,
the automorphism of the free group on $n$ free generators. If $\rho$ is a
complex linear representation of $\mathrm{Aut}(F_n)$ then $\rho(\mathrm{Inn}(F_n))$
is virtually soluble, where $\mathrm{Inn}(F_n)$ is the subgroup of inner
automorphism of $F_n$.
\end{conjecture}

\section{Theorem \ref{nori}}

\label{nor}

Our sketch of the proof of Theorem \ref{nori} follows the argument in
\cite{SG}, Window 9.
\medskip

Suppose that $\Gamma \leq G_{\mathbb Z_S}$ is Zariski dense in the simply connected $\mathbb Q$-simple
algebraic group $G$. Now $G$ may not be absolutely simple, but in any case
there is a number field $k$ and an absolutely simple group $H$ defined over
$k$ such that $G= \mathcal R_{k/\mathbb Q}(H)$.
We have $G_\mathbb Q= H_k$, and for each prime $p$
\[ G_{\mathbb Z_p}= \prod_{j} H_{\mathcal O_{\mathbf p_j}} \]
where $p\mathcal O= \prod_{j} \mathfrak p_j^{e_j}$ is the factorization of the principal
ideal $(p)$ in $\mathcal O$. This means that $k \otimes \mathbb Q_p= \prod_j
k_{\mathfrak p_j}$. 

From now on assume that the prime $p$ is unramified in $k$, i.e. all $e_j=1$. In addition assume that $G$ has good reduction mod $p$. This holds for all but finitely many rational primes $p$ (see Proposition \ref{redp}). 

Since $L(G)$ is $\mathbb Q$-defined we have that $L(G)_{\mathbb Q_p} = L(G)
\otimes \mathbb Q_p$. 
Therefore $L(G)_{\mathbb Q_p}= \prod_j L(H)_{k_{\mathbf p_j}}$. Similarly since $p$ is unramified

\begin{equation} \label{deco} 
L(G)_{\mathbb F_p} = \prod_j  L(H)_{\mathcal O/\mathbf p_j} \ \mathrm{\
 and} \end{equation}
 
\[ G_{\mathbb F_p}= \prod_j H_{\mathcal O/\mathbf p_j}.\]
\medskip

The group $H$ is absolutely simple so for almost all primes $p$ the Lie algebras
$L(H)_{\mathcal O /\mathbf p_j}$ are simple and the groups $H_{\mathcal O/\mathbf p_j}$ are quasisimple.
\bigskip

\textbf{Step 1}: Let $D_p$ be the closure of $\Delta$ in the $p$-adic analytic group
$G_{\mathbb Q_p}$. Since $\Delta$ 
is Zariski-dense in $G$ then by Proposition \ref{ideals} the Lie algebra
of $D$ is an ideal of the Lie algebra $L(G)_{\mathbb{Q}_p}$ 
of $G_{\mathbb Q_p}$. But $\Delta \leq G_\mathbb Q$, so the Lie algebra $L(D_p)$
is defined over $\mathbb Q$. Hence the projections of
$L(D_p)$ in each of the factors $L(H)_{k_{\mathbf p_j}}$ of $L(G)_{\mathbb{Q}_p}$
are isomorphic. So for almost all primes $p$ we have $L(D_p)=L(G)_{\mathbb Q_p}$ which means that  $D_p$ is an open subgroup of $G_{\mathbb
Q_p}$ for almost every $p$ (see Proposition \ref{openlie}). In fact, since we are assuming $p \not \in S$,
we have $\Delta \subset G_{\mathbb Z_p}$, and so $D_p$ is an open subgroup of
the compact open subgroup $G_{\mathbb Z_p}$.

Next we want to prove that for almost all primes $p$ our group $\Delta$ is dense in $G_{\mathbb
Z_p}$.

\medskip

\textbf{Step 2}: For almost all primes the Frattini subgroup of $G_{\mathbb Z_p}$ is contained in the kernel of $G_{\mathbb{Z}_p} \rightarrow G_{\mathbb F_p}$. In follows that a subgroup $\Delta$ is dense in $G_{\mathbb Z_p}$ if and only if $\Delta$ maps onto $G_{\mathbb F_p}$.
This is proved in \cite{SG},
Window 9, Proposition 7 using the structure of the finite
images of the $p$-adic analytic group $G_{\mathbb Z_p}$.

\medskip

\textbf{Step 3}: We shall prove that $D_p=G_{\mathbb Z_p}$ for almost all primes $p$.
By Step 2 it is enough to show that $\Delta$ maps onto $G_{\mathbb F_p}$ for almost all primes $p$. 

Let $\pi_p$ be the projection of $G_{\mathbb Z_p}$ onto $G_{\mathbb F_p}$
and further
let $\pi_j$ and and $\tau_j$ be the projections of $G_{\mathbb Z_p}$ and $L(G)_{\mathbb Z_p}$ onto their direct factors $H_{\mathcal O/\mathbf p_j}$ and $L(H)_{\mathcal
O/\mathbf p_j}$ respectively.

\bigskip

At this stage we need the following

\begin{proposition}
\label{key} Let $\Gamma$ be a subgroup of $G_{\mathbb{F}_{p}}$ such that

(a) For all $j$ the image $\pi_{j}(X)$ of $\Gamma$ in $H_{\mathcal{O}%
/\mathfrak{p}_{j}}$ has order divisible by $p$, and

(b) Every subspace of $L(G)_{\mathbb{F}_{p}}$ invariant under $\Gamma$ is an ideal.

Then provided $p$ is sufficiently large compared to $\dim G$ we have 
$\Gamma= G_{\mathbb{F}_{p}}$.
\end{proposition}

Let us check that the conditions (a) and (b) above are satisfied for the group
$\pi_{p}(\Delta) \leq G_{\mathbb{F}_{p}}$, for almost all primes $p$. \medskip

Suppose that (a) fails for a set $A$ of infinitely many primes. Then there is
$j=j_{p}$ such that $\pi_{j_{p}}(\Delta)$ has order coprime to $p$ and so is a
completely reducible subgroup of $\mathrm{GL}_{n}(\mathbb{F}_{p})$, where $n$
depends only on $G$ and not on $p$. A variation of Jordan's Theorem \cite{jordan} then says that there
is a number $f=f(n)$ such that $\pi_{j,p}(\Gamma)$ has an abelian subgroup of
index at most $f$.

Since the set $A$ of rational primes is infinite we have
\[
G_{\mathbb{Z}_{S}} \cap\bigcap_{p \in A} \ker\pi_{j_{p}} = \{1\}
\]

This implies that $\Delta$ itself is virtually abelian (it is finitely generated
so it has only finitely many subgroups of index at most $f(n)$). But $\Delta$
is Zariski-dense in the $\mathbb{Q}$-simple algebraic group $G$: contradiction.

So condition (a) of Proposition \ref{key} holds for almost all primes.

Condition (b) is immediate: $H$ is absolutely simple and so for almost all
primes each of the $L(H)_{\mathcal{O}/\mathfrak{p}_{j}}$ is a simple module
for $H_{\mathcal{O}_{/}\mathfrak{p}_{j}}$. Since $\Delta$ is Zariski- dense in
$H_{k}$ the group $\mathrm{Ad}(\Delta)$ spans $\mathrm{End}_{k}L(H)_{k}$ so for almost all
primes $\mathrm{Ad}(\pi_{j}(\Delta))$ spans End$_{\mathcal{O}/\mathfrak{p}_{j}%
}L(H)_{\mathcal{O}/\mathfrak{p}_{j}}$. This means that each summand
$L(H)_{\mathcal{O}/\mathfrak{p}_{j}}$ of $L(G)_{\mathbb{F}_{p}}$ is a simple
module for $\pi_{p}(\Delta)$. So the decomposition of
$L(G)_{\mathbb{F}_{p}}$ into minimal Lie ideals is also a decomposition into
irreducible $\mathbb{F}_{p}\pi_{p}(\Delta)$-modules. So every irreducible
module for $\pi_{p}(\Delta)$ in $L(G)_{\mathbb{F}_{p}}$ is an ideal, proving
that (b) holds. \medskip

\textbf{Step 4} We now know that the closure $\overline{\Delta}$ of $\Delta$
in $G_{\widehat{\mathbb{Z}}_{S}}= \prod_{ p \not \in S} G_{\mathbb{Z}_{p}}$
projects onto all but finitely many of the factors $G_{\mathbb{Z}_{p}}$. Now
it is easy to show (see Exercise 16) that in this case $\overline{\Delta}$
contains their Cartesian product. Combined with Step 1 (which says that
$\overline{\Delta}$ projects onto an open subgroup in each of the remaining
factors) we easily see that $\overline{\Delta}$ is open in $G_{\widehat
{\mathbb{Z}}_{S}}$.

\subsection{Proposition \ref{key}}

There are now at least three different proofs of Proposition \ref{key}. One is
by Matthews, Vaserstein and Weisfeiler \cite{matt}, it uses the Classification
of the Finite simple groups to deduce properties of a proper subgroup of
$G_{\mathbb{F}_{p}}\leq\mathrm{GL}_{n}(\mathbb{F}_{p})$ which are
incompatible with (a) and (b).

There is also a proof using logic by Hrushovkii and Pillay \cite{Hrush}.

We shall focus on the original proof by Nori \cite{nor}. It studies
unipotently generated algebraic groups and their Lie algebras in large finite
characteristic $p$. This is motivated by the construction of 
the Chevalley groups described on section \ref{chev}. Recall that the adjoint chevalley group $\overline G$
is generated by certain automorphisms $\exp (\mathrm{ad}(x))$ for 
certain $\mathrm{ad}$-nilpotent elements $x$ of the Lie algebra of $G$. If we fix such an element $x$ 
then the set 
\[ \{ \exp(\mathrm{ad}(tx)) \ | \ t \in K  \} \]
is a unipotent subgroup of $\overline{G}$ and is isomorphic to $\mathbb G_+$.

Nori generalizes this situation in two directions: He proves an analogue of this 
not just for algebric groups but for Zariski-dense subgroups of 
$\mathrm{GL}_n(\overline{\mathbb F_p})$ and secondly, he does this not 
just in the algebraic closure $\overline{\mathbb F_p}$ of $\mathbb F_p$ but in the finite field $\mathbb F_p$ 
(provided $p$ is large enough compared to $n$). \medskip   

The details are as follows: \medskip

For a group $\Gamma\leq\mathrm{GL}_{n}(\mathbb{F}_{p})$ let $\Gamma^{+}$ be
the subgroup generated by its unipotent elements. When $p\geq n$ these are
just the elements of order $p$ in $\Gamma$. Similarly for an algebraic group
$G\leq\mathrm{GL}_{n}(K)$ let $G^{+}$ be the subgroup generated by its
unipotent elements.


Now for an element $g\in\mathrm{GL}_{n}(\mathbb{F}_{p})$ of order $p$ let
$X_{g}$ be the unipotent 1-dimensional algebraic group over $\overline
{\mathbb{F}}_{p}$ generated by $g$. In other words define
\[
X_{g}=\left\{  g^{t}:=\sum_{i=0}^{p}\left(  \overset{t}{i}\right)
(g-1)^{i}\quad|\quad t\in\overline{\mathbb{F}}_{p}\right\}  ,
\]
where $\overline{\mathbb{F}}_{p}$ is the algebraic closure of $\mathbb{F}_{p}%
$. Note that $X_{g}$ is defined over $\mathbb{F}_{p}$ and is isomorphic to the
additive group of the field $\overline{\mathbb{F}}_{p}$.

Now, given $\Gamma\leq\mathrm{GL}_{n}(\mathbb{F}_{p})$ define the algebraic
group $T=T(\Gamma)$ as
\[
T=\langle X_{g}\ |\ \forall g\in\Gamma,g^{p}=1\rangle\leq\mathrm{GL}%
_{n}(\overline{\mathbb{F}}_{p}).
\]
Recall that the subgroup generated by a collection of closed connected
subgroups is closed and connected, so $T$ is indeed a connected algebraic group.
Observe that since $X_g$ is the smallest connected algebraic group containing $g$ and $g \in G_{\mathbb F_p}$ it follows that $X_g \leq G$ and hence $T \leq G$.

Nori's main result is that in the above setting we have
\[
\Gamma^{+} = (T_{\mathbb{F}_{p}})^{+}%
\]
provided $p$ is large enough compared to $n$.

Now, it is known that for large primes $p$ one has 
\[ (T_{\mathbb{F}_{p}})^{+}=(T_{\mathbb{F}_{p}}).\] 

So $\Gamma^+$ is the group of $\mathbb{F}_{p}$-rational
points of the connected algebraic group $T$. 

Now, suppose that condition (b)
of Proposition \ref{key} holds. Clearly $\Gamma$ normalizes the algebraic group $T \leq G$ since 
$(X_g)^{\gamma}= X_{g^\gamma}$ for any $\gamma, g \in \Gamma$ with $g^p=1$. Therefore the Lie algebra 
$L(T) \leq L(G)$ of $T$ is normalized by $\Gamma$.

 It follows that the subspace $L(T)_{\mathbb{F}_{p}}$
of $L(G)_{\mathbb{F}_p}$ is invariant under $\Gamma$ and so it is an ideal of $L(G)_{\mathbb{F}%
_{p}}$. Not only that, $L(T)$ is defined over $\mathbb{F}_{p}$ and so its
projections on the direct factors of $L(G)$ are isomorphic. In the same way as
in Step 1 above we deduce that $L(T)=L(G)$ and since both $G$ and $T \leq G$
are connected we have $T=G$. So%

\[
\Gamma\geq\Gamma^{+}= T_{\mathbb{F}_{p}}=G_{\mathbb{F}_{p}} \geq\Gamma
\]
giving that $\Gamma= G_{\mathbb{F}_{p}}$ as required.

\section{Exercises}

1. Show that every open set in $K^{n}$ can be regarded as closed affine set in
some $K^{m}$, $m\geq n$. \medskip

2. Prove that $\dim V$ for an irreducible affine variety $V$ is the largest
$d$ such that we can find a chain $\emptyset\not = V_{1} \subset V_{2}
\subset\cdots V_{d} \subset V$ of distinct irreducible closed subvarieties
$V_{i}$ in $V $. You may use any of the equivalent definitions of $\dim V$ in
\S 2.1. \medskip

3. (Proposition \ref{compact}) \label{exerc} Show that each affine variety is a compact topological space and that in
fact it satisfies the descending chain condition on closed subsets. \medskip

A subset $X \subset V$ of an affine variety $V$ is \emph{constructible} if it
can be obtained from the open or closed subsets of $V$ by a finite process of forming 
unions and intersections. A theorem of Chevalley says that an image of a
constructible set under a morphism of varieties is constructible. \medskip

4. (\cite{Wef}, Lemma 14.10.) Prove that a constructible (abstract) subgroup $H$ of a linear algebraic
group $G$ is in fact closed and so is algebraic. Deduce with Chevalley's theorem
that an image of an algebraic group under a homomorphism is an algebraic
group. 

5. (\cite{Wef}, Lemma 14.14) Let $G$ be a linear algebraic group and $(X_{i})_{i \in I}$ be a family of
constructible irreducible subsets of $G$ each containing the identity. Show
that $X_{i}$ together generate a closed irreducible subgroup of $G$. Hence
deduce that if $G$ is connected, then the derived subgroup $G^{\prime}= \langle[x,y] | \ x,y \in G
\rangle$ is both closed and connected. \medskip

6. Suppose that $k/k_{0}$ is a finite extension of fields and $H=\mathcal{R}%
_{k/k_{0}}(G)$. Show that $H$ is $K$-isomorphic to
\[
G^{\sigma_{1}} \times G^{\sigma_{2}} \times\cdots\times G^{\sigma_{d}}%
\]
where $\sigma_{i}$ are all the embeddings of $k$ in $K$ which fix the elements of $k_{0}$ and
$G^{\sigma_{i}}$ is the algebraic group defined by the ideal $I^{\sigma_{i}}$
where the ideal $I$ defines $G=V(I)$ as a variety in $M_{n}(K)$. Hint: 
use the map $\lambda$ on page \pageref{chin} and the isomorphism (\ref{chin}). \medskip

7. \label{ex7} Let $G$ be the multiplicative group of norm one quaternions
defined over $\mathbb{Q}$. For example we can take $G$ in its left regular representation

\[
G=\left\{
\begin{pmatrix}
a & - b & -c & -d \\ b & a & d & -c \\ c & -d & a & b \\
d & c & -b & a 
\end{pmatrix}
\ \left| {}\right.  \ a^{2}+b^{2}+c^{2}+d^{2}=1\right\}
\]

Show that $G$ is $\mathbb Q(i)$-isomorphic to $\mathrm{SL}_2$ but it is not
$\mathbb Q$-isomorphic to it. Hint: Send the $4 \times  4$ matrix with first
column $a,b,c,d$ as above to

\[ \begin{pmatrix} a+ib & -c+id \\ c +id & a-ib \end{pmatrix}. \]

8. Show that if $G=\mathrm{SL}_{n}(K)$ then $L(G)=\mathrm{sl}_{n}(K)$, the Lie algebra of
matrices of trace $0$ in $M_{n}(K)$.\medskip

9. Show that $\Gamma=\mathrm{SL}_{2}(\mathbb{Z})$ does not have the
generalized congruence subgroup property. You may use that $\Gamma$ has a
nonabelian free subgroup of finite index. \medskip

10. Show that $\mathrm{SL}_{n}(\mathbb{Z})$ has the strong approximation property.
(Hint: use the fact that for a finite ring $R$ the group $\mathrm{SL}_{n}(R)$ is generated by
elementary matrices.) \medskip

11. Show that $\mathrm{PGL}_{2}(\mathbb{Z})$ fails to have the strong
approximation property (as an arithmetic subgroup of $G=\mathrm{PGL}_{2}$).
\medskip

12. Show that if a connected linear algebraic group $G$ is not soluble then it
maps onto a simple algebraic group. (Hint: Let $M=\mathrm{Rad }\ G$ be the
soluble radical of $G$. Then $G/M$ is semisimple.) \medskip

13. Suppose that $\Gamma$ is a Zariski-dense subgroup of a connected algebraic
group $G$ and that $\Delta$ is a subgroup of finite index in $\Gamma$. Show
that $\Delta$ is also Zariski-dense in $G$. \medskip

14. Suppose that $G\leq\mathrm{GL}_{n}(K)$ is a connected algebraic group
which has a normal subgroup $N$ which preserves a one-dimensional subspace
$\langle v\rangle$. Show that either $N$ acts as scalars or else $G$
stabilizes a nontrivial subspace of $K^{n}$. \medskip

15. Show that a connected soluble algebraic group $G\leq \mathrm{GL}_{n}(K)$
has a common eigenvector. Deduce that $G$ is triangularizable and hence prove
Theorem \ref{lk}. (Hint: use Exercise 14 with $G^{\prime}$ in place of $N$.)
\medskip

16. Suppose that $L$ is a closed subgroup of $K=\prod_{ p\in A }
G_{\mathbb{Z}_{p}}$ for some set $A$ of primes, where $G$ is a $\mathbb{Q}%
$-simple connected and simply connected algebraic group.

(a) Show that if $p$ is sufficiently large then if $L$ maps onto the direct
factor $G_{\mathbb{Z}_{p}}$ of $K$ then in fact it contains it.

(b) On the other hand if $A$ is finite set of primes and $L$ maps onto an open
subgroup of each factor $G_{\mathbb{Z}_{p}}$ of $K$ show that then $L$ is an
open subgroup of $K$. \medskip

17. Show that for any algebraic group $G$ in characteristic 0 the group
$G^{+}$ generated by its unipotent elements is connected. (Hint: use exercise 5) \medskip

18. Using Theorems \ref{reductive} and \ref{K-simple} show that if a connected
algebraic group consists of semisimple elements then it is a torus. 
(Hint: a nontrivial semisimple group contains a copy of $\mathrm{SL}_2$ or 
$\mathrm{PSL}_2$.)\medskip

19. (\cite{one4all}) Let $n>1$ be an integer. Show using Strong Approximation that there is a 
finite set $A$ of rational primes with the following property: If $S \subseteq \mathrm{SL}_n(\mathbb Z)$ is 
a subset whose image generates $\mathrm{SL}_n(\mathbb F_p)$ for some prime $p \not\in A$, then for almost all 
primes $q$, the image of $S$ in $\mathrm{SL}_n(\mathbb F_q)$ generates $\mathrm{SL}_n(\mathbb F_q)$. 
Generalize this to any absolutely simple, connected, simply connected group $G$ defined over $\mathbb Z$.

\end{document}